\documentclass[11pt,reqno]{amsart}
\usepackage[T1]{fontenc}
\usepackage[utf8]{inputenc}
\usepackage{lmodern}
\usepackage{amsmath,amssymb,amsthm}
\usepackage{array}
\usepackage{placeins}
\usepackage[margin=1.03in]{geometry}
\usepackage{microtype}
\usepackage{tikz}
\usetikzlibrary{arrows.meta,positioning,decorations.pathreplacing,calc}
\definecolor{leftbranch}{RGB}{36,92,148}
\definecolor{rightbranch}{RGB}{163,82,23}
\usepackage[hyphens]{url}
\usepackage[hidelinks]{hyperref}
\hypersetup{pdftitle={The geometric case of the Erd\H{o}s similarity conjecture: A human exposition of the OpenAI proof},pdfauthor={Alex Iosevich and Alexia Yavicoli}}
\newenvironment{gridtable}{\begin{center}\vspace{6pt}\small}{\par\vspace{6pt}\end{center}}
\newtheorem{theorem}{Theorem}[section]
\newtheorem{proposition}[theorem]{Proposition}
\newtheorem{lemma}[theorem]{Lemma}
\theoremstyle{remark}

\numberwithin{equation}{section}
\allowdisplaybreaks[2]
\title[The geometric case of the Erd\H{o}s similarity conjecture]
{The geometric case of the Erd\H{o}s similarity conjecture\\A human exposition of the OpenAI proof}

\author{A. Iosevich}
\address{Department of Mathematics, University of Rochester, Rochester, NY, USA}
\email{iosevich@gmail.com}

\author{A. Yavicoli}
\address{Department of Mathematics, University of British Columbia, Vancouver, BC, Canada}
\email{yavicoli@math.ubc.ca}

\thanks{The theorem and the finite tree construction are due to OpenAI, in the manuscript dated October 5, 2026 and released on October 6, 2026. This exposition expands the proofs and compares their ingredients with earlier work.
A.~I. was supported in part by National Science Foundation grant
DMS-2154232. A.~Y. was supported in part by the Natural Sciences and Engineering Research
Council of Canada (NSERC), grants GR030571 and GR030540.}
\begin{document}
\begin{abstract}
For a fixed $q\in(0,1)$, can a set of positive Lebesgue measure avoid every translated and dilated copy of $\{q^n:n\ge1\}$? OpenAI's release of October 6, 2026 shows that it can, with a compact subset of $[0,1]$ whose measure is arbitrarily close to one. We give a detailed exposition of the proof, using only finite probability and basic measure theory. The proof removes a small open set intersecting every copy, obtained from a random construction: fair random choices, attached to the cells of nested grids, route each point down a finite tree, and at the leaf reached a further random choice removes the point with probability $p$. The removed set therefore has expected density $p$, however large the tree, while the branching gives each copy many independent chances to intersect it. Each chance depends on only part of the tree, so its outcome changes at only a controlled number of dilations, and a union bound over finitely many representative dilations then covers every dilation.
We also relate the proof to the earlier random-cell arguments going back to Kolountzakis (Bull. London Math. Soc. 29 (1997), 415--424) and to recent work on two-fold sumsets (arXiv:2607.03584) and on Rajchman measures (arXiv:2609.04456), and we show how the method of conditional expectations makes the construction deterministic.
\end{abstract}
\maketitle
\tableofcontents
\clearpage

\section{The question and the conclusion}\label{sec:question}

Suppose that a measurable subset of the real line has positive Lebesgue measure. Must it contain a translated and dilated copy of a prescribed set?
It is easy to prove that the answer is positive if the prescribed set is finite.
When the prescribed set is infinite, the answer is believed to be negative. The problem is of particular interest for sequences. Any fixed finite part of a sequence is easy to find once the dilation has been made small enough, but keeping the same translation and dilation for the entire sequence is a different matter.

Write $m$ for Lebesgue measure. A set $A\subset{\mathbb R}$ is called \emph{measure universal} if every measurable set of positive measure contains
\[
x+sA=\{x+sa:a\in A\}
\]
for some $x\in{\mathbb R}$ and some real $s\ne0$, negative values of $s$ allowing reflection. The Erd\H{o}s similarity conjecture~\cite[problem~4.33.7*]{Erdos1974} asserts that no infinite set is measure universal, and the geometric case concerns
\[
G_q=\{q^n:n\ge1\},
\]
where $q\in(0,1)$ is fixed. Here and below, ${\mathbb N}=\{1,2,3,\ldots\}$, and the indices $n$ are natural numbers. For $x\in{\mathbb R}$ and real $s\ne0$, we call $x$ the \emph{base point} of the copy $x+sG_q$. It is the limit of the points $x+sq^n$ as $n\longrightarrow\infty$, and it does not belong to the copy.

\begin{theorem}[OpenAI \cite{OpenAI2026}]\label{thm:main}
Let $q\in(0,1)$ and $\eta\in(0,1)$. There is a compact set $E\subset[0,1]$ with $m(E)>1-\eta$ such that, for every $x\in{\mathbb R}$ and every real $s\ne0$, there is an $n\ge1$ for which $x+sq^n\notin E$.
\end{theorem}

The set $E$ may depend on $q$ and $\eta$, so that the theorem does not assert that a single set avoids all geometric ratios simultaneously. It proves the geometric case of the conjecture, and not the conjecture for arbitrary infinite sets.
For the history of the problem, partial solutions and its variants, see the papers by Svetic \cite{Svetic2001}, and by Jung, Lai, and Mooroogen~\cite{JungLaiMooroogen2025}.

Let us first justify the finite observation with which we began. Suppose $F=\{a_1,\ldots,a_N\}$ is finite and $A$ is measurable with $0<m(A)<\infty$. Set
\[
A_s=\bigcap_{j=1}^N(A-sa_j).
\]
Continuity of translations in $L^1({\mathbb R})$ means that
\[
m(A\mathbin{\triangle}(A-h))\longrightarrow0
\]
as $h\longrightarrow0$, where $\triangle$ denotes symmetric difference. Since
\[
m(A\setminus A_s)\le\sum_{j=1}^Nm(A\setminus(A-sa_j)),
\]
we have $m(A_s)>0$ for all sufficiently small $|s|$, and any $x\in A_s$ satisfies $x+sF\subset A$. If $m(A)=\infty$, apply the argument to a bounded subset of $A$ of positive measure. Finite patterns are thus universal. Since the size of the permissible dilation can shrink as more points are added, this argument says nothing about a single copy of an infinite sequence.

\subsection{What the reader will need}

The proof of Theorem \ref{thm:main} uses finite collections of independent random variables, integration over a unit interval, compactness, and approximation of a measurable set from outside by open sets. A Bernoulli-$p$ random variable equals one with probability $p$ and zero with probability $1-p$, and we write ${\mathbb P}$ and ${\mathbb E}$ for probability and expectation. Conditioning on specified table entries will mean fixing those entries and averaging over all the entries which remain. No infinite random process is used.

The main proof occupies sections~\ref{sec:periodic}--\ref{sec:finish}. We first identify the small periodic set which would imply the theorem, then explain why tests on a single grid, checking whether points of a copy lie in randomly selected cells, are insufficient, and construct a tree which avoids this difficulty, the geometry of its grids being established before the probability calculation. We then pass from one dilation to all dilations in a bounded interval, use a descent calculation to assemble the stopping estimates, and repair a small set of exceptional base points. Section~\ref{sec:descent} explains the connection with Burstein--Iosevich--Krause, and section~\ref{sec:comparison} returns to the two-fold sumset and Rajchman papers; that comparison can be read independently after the proof, and Fourier analysis appears only there.

\subsection{The relation with earlier work}

Falconer~\cite{Falconer1984} and Eigen~\cite{Eigen1985} proved nonuniversality for positive decreasing sequences $(a_n)_{n\in{\mathbb N}}$ tending to zero with $\frac{a_{n+1}}{a_n}\longrightarrow1$ as $n\longrightarrow\infty$, Humke and Laczkovich~\cite{HumkeLaczkovich1998} reformulated the problem in terms of finite coverings and obtained further conditions on the gaps between successive terms, and Bourgain~\cite{Bourgain1987} approached the problem through additive structure, proving that the sum of three infinite sets is never measure universal. The two-fold sumset results discussed below belong to this line of work, but their proofs follow a different approach from Bourgain’s.

The random-cell argument used here follows the work of Kolountzakis~\cite{Kolountzakis1997}. Chleb\'ik~\cite[section~5]{Chlebik2015} and Kolountzakis--Papageorgiou~\cite[section~3.1]{KolountzakisPapageorgiou2025} count the finitely many ways a dilation can change cell membership and then integrate the probability that a base point is exceptional.

The dyadic sequence was already studied by random cells in the reports of David Solymosi~\cite[section~3]{Solymosi2011} and Foster Tom~\cite{Tom2015}, both describing undergraduate research with Malabika Pramanik. Solymosi studies finite choices of the dilation and the measure of the remaining exceptional base points, the two quantities controlled in sections~\ref{sec:scales} and~\ref{sec:finish} below. Tom considers nested dyadic intervals, of the kind used for the grids of section~\ref{sec:geometry}, and asks for checks which work independently of the dilation. Neither report gives a complete proof of the dyadic case, but both identify the difficulty overcome below, and Tom's lemma~0.2 is the principle used in section~\ref{sec:finish} to repair exceptional base points: an open neighborhood of these base points removes the remaining copies of a null sequence.

Two of the papers cited by OpenAI are particularly close to the questions we ask below. The two-fold sumset paper of Mora Cu\'ellar, Iosevich, Kulkarni, Rojas Aravena, and Yavicoli~\cite{IKMRY2026Sumset} obtains enough separated points to overcome the number of cell configurations met as the dilation varies in a random construction, while the Rajchman paper of Iosevich, Kulkarni, Mora Cu\'ellar, Rojas Aravena, and Yavicoli~\cite{IKMRY2026Rajchman} constructs small periodic sets which intersect every translate at the required dilations and then combines them with a summable measure bound. We shall identify these common features at the places where they arise, the finite tree and its estimate for the number of relevant cell configurations being the additional mechanism in the OpenAI proof.

There are related results in which either the size of the avoiding set or the allowed maps are different. Cruz, Lai, and Pramanik~\cite{CruzLaiPramanik2023} construct sets of Hausdorff dimension one and Lebesgue measure zero which avoid affine copies of certain null sequences together with their limit point, including $G_q\cup\{0\}$, while Feng, Lai, and Xiong~\cite{FengLaiXiong2024} show that every geometric null sequence can be embedded into every measurable set of positive measure by a bi-Lipschitz map, that is, one which changes distances by at most fixed multiplicative factors. Both positive Lebesgue measure and the restriction to affine copies therefore matter in the question considered here.

\section{Begin with the set we want to remove}\label{sec:periodic}

To build a large set avoiding every copy of $G_q$, it is enough to remove a small open set which intersects every copy, and we call such a set a \emph{hitting set}. This change of viewpoint is useful because different families of copies can be handled by different small sets, whose measures can then be added.

For a measurable set $H$ satisfying $H+1=H$, define its \emph{density} by
\[
\rho(H)=m(H\cap[0,1)).
\]
The equality $H+1=H$ means that $H$ is $1$-periodic. We will prove the following statement, keeping $q$ fixed throughout.

\begin{proposition}\label{prop:periodic}
For every $p\in(0,1)$ there is an open $1$-periodic set $H\subset{\mathbb R}$ such that $\rho(H)\le6p$ and
\begin{equation}\label{eq:hitting}
\text{for every }x\in{\mathbb R}\text{ and }t\in[1,2],
\text{ some }n\ge1\text{ satisfies }x+tq^n\in H.
\end{equation}
\end{proposition}

The parameter $p$ will control the amount removed, while the interval $[1,2]$ is only a normalization of the dilation parameter and has no arithmetic relation to $q$.

\begin{proof}[Proof of theorem~\ref{thm:main} from proposition~\ref{prop:periodic}]
For each integer $k$, choose a set $H_k$ from the proposition with
\[
p_k=\frac{\eta}{64}4^{-|k|}.
\]
Define
\[
C=\bigcup_{k\in{\mathbb Z}}(2^kH_k\cup(-2^kH_k)),
\]
and put $E=[0,1]\setminus C$. Since $C$ is open, $E$ is compact.

The period of either signed dilation $\pm2^kH_k$ is $2^k$, and its measure in one period is $2^k\rho(H_k)$. If $k\le0$, then the unit interval is a union of exactly $2^{-k}$ periods, while if $k>0$, it is contained in the period interval $[0,2^k]$. Consequently,
\[
m((\pm2^kH_k)\cap[0,1])\le\max\{1,2^k\}\rho(H_k).
\]
Countable subadditivity now yields
\begin{equation}\label{eq:global-budget}
\begin{aligned}
m(C\cap[0,1])
&\le12\sum_{k\in{\mathbb Z}}p_k\max\{1,2^k\}\\
&=\frac{12\eta}{64}\Bigl(1+\sum_{j=1}^{\infty}4^{-j}+\sum_{k=1}^{\infty}2^{-k}\Bigr)
=\frac{7\eta}{16}<\eta.
\end{aligned}
\end{equation}
Thus $m(E)>1-\eta$.

Given $s>0$, choose $k\in{\mathbb Z}$ and $t\in[1,2)$ with $s=2^kt$. The rescaled copy $2^{-k}(x+sG_q)=2^{-k}x+tG_q$ has base point $2^{-k}x$, and applying~\eqref{eq:hitting} to $H_k$ at this base point, we find $n\ge1$ with
\[
2^{-k}x+tq^n\in H_k.
\]
Multiplication by $2^k$ shows that $x+sq^n\in C$. If $s<0$, apply the positive case to $-x$ and $-s$ and then reflect; since the set $C$ is symmetric, this also produces a point of $x+sG_q$ in $C$.
\end{proof}

This is the first direct connection with the Rajchman paper~\cite[proposition~3.2]{IKMRY2026Rajchman}, both arguments reserving summable amounts of measure for open periodic sets which handle specified dilations. What remains to be proved is the local hitting statement. The Rajchman argument obtains it from approximate density modulo one, whereas here we must produce it for a countable geometric sequence by another method.

\section{Why one grid does not settle the question}\label{sec:obstacle}

Imagine dividing the line into small periodic cells and selecting each cell independently with probability $p$. We use the word \emph{test} for a check attached to a chosen point of a copy, at a given dilation, with two possible outcomes, success or failure; a successful test certifies that the point belongs to the random set, so a copy intersects the random set as soon as one of its tests succeeds. In the single-grid construction considered here, the test simply checks whether the cell containing the point has been selected. If $N$ points of a copy occupy different cells, then the probability that all of them miss the selected cells is $(1-p)^N$, and the elementary inequality $1-u\le e^{-u}$ implies
\[
(1-p)^N\le e^{-pN}.
\]
For a fixed translation and dilation, increasing $N$ therefore looks promising.

The difficulty is that the dilation is not fixed. Fix the base point $x$ and let the dilation $t$ vary over $[1,2]$. Each of the chosen points $x+tq^n$ moves along an interval, changing its cell when it crosses a grid boundary. The \emph{cell configuration} of the chosen points is the ordered list of the cells they occupy, each cell being recorded by its position modulo one; from section~\ref{sec:geometry} on, this is the ordered list of their cell labels. Only finitely many cell configurations occur, and a bound valid for every dilation must account for all of them. If there are $\Lambda$ cell configurations, and if at every dilation under consideration the $N$ tested points occupy distinct cells, the union bound estimates the probability of a miss at some dilation by $\Lambda e^{-pN}$. The union bound is simply the inequality saying that the probability of a union of events is at most the sum of their probabilities, and it requires no independence. We therefore need enough tests to outweigh the growth of $\Lambda$ while keeping $p$ small.

For a geometric sequence, the exponential decay $e^{-pN}$ of the failure probability competes with the number $\Lambda$ of cell configurations. Once the cells are fine enough to separate the chosen points, $\Lambda$ grows at least exponentially in $N$, because the smallest gap between $N$ powers of $q$ is exponentially small compared with their spread. A precise way to see the obstruction is to take any $N\ge2$ distinct powers
\[
a_1>a_2>\cdots>a_N
\]
of $q$, and define their normalized minimum gap by
\[
\delta=\frac{\min_{1\le j<N}(a_j-a_{j+1})}{a_1-a_N}.
\]
Since $a_{N-1}\le q^{N-2}a_1$ and $a_1-a_N\ge(1-q)a_1$, we have
\begin{equation}\label{eq:gap-obstruction}
\delta\le\frac{q^{N-2}}{1-q}.
\end{equation}
In particular $-\log\delta$ cannot be $o(N)$, where $u_N=o(N)$ means that $\frac{u_N}{N}\longrightarrow0$ as $N\longrightarrow\infty$, and $\log$ denotes the natural logarithm.

This is exactly the obstruction discussed in~\cite[remark~15]{IKMRY2026Sumset}, where large finite subsets are used for which the logarithm of the reciprocal gap is small compared with the number of points, the additional summand there supplying many separated cross-sums. A single geometric sequence has no such supply, and the inequality~\eqref{eq:gap-obstruction} also explains why passing to more widely spaced powers does not repair that particular criterion.

The OpenAI construction overcomes this obstruction with a finite rooted tree which organizes the random choices. Each leaf carries its own random selection of cells, each cell being selected with probability $p$. Separate random choices, independent of these selections, send each real point along a single path from the root to exactly one leaf, and the point belongs to the random set if its cell is selected at that leaf. Consequently, each point belongs to the random set with probability $p$, however many branches the tree has, and the expected density of the set remains $p$.

The branches also organize the points used to detect a copy. The finite set of sequence indices used in the construction is divided into intervals of consecutive integers, and one such index interval is assigned to each edge of the tree. For a copy $x+tG_q$, an index interval supplies the points $x+tq^n$ whose indices lie in it. The route of the base point $x$ determines which index intervals enter the probability estimate: those on the edges leaving one particular vertex of that route. For all base points outside a set of small density, the points of the copy coming from these intervals lie so close to $x$ that they make the same random choices as $x$ up to that vertex, while from there on they read their own random choices, independent of those of $x$. For these points we construct local tests, each of which can certify membership in the random set by inspecting only the choices on the corresponding edge and in the subtree below it. In this way the branching supplies many independent opportunities for a successful test.

The remaining issue is to handle every dilation $t\in[1,2]$. As $t$ varies, a local test can change its outcome only when its point crosses a boundary of one of the grids that the test consults, so grids elsewhere in the tree do not enter the count of cell configurations for that test. The lengths of the index intervals are chosen to keep the number of these crossings under control. This is the advantage of the construction: increasing the branching supplies more opportunities to intersect the random set, while the arrangement of the index intervals keeps the number of relevant cell configurations small enough for the union bound to succeed. Sections~\ref{sec:geometry}--\ref{sec:scales} establish these geometric and probabilistic properties.

\section{A tree, its index intervals, and its grids}\label{sec:geometry}

We now describe the deterministic part of the construction, choosing the numerical parameters only after all the estimates are available, so as to avoid a circular choice of constants.

Fix natural numbers $M\ge2$, $d\ge1$, $g\ge1$ and $r_0\ge1$. The \emph{depth} of a vertex in a rooted tree is the number of edges on its path from the root. We use the finite rooted tree with levels $0,\ldots,d$ in which the root has depth $0$, every vertex of depth less than $d$ has exactly $M$ children, numbered $1,\ldots,M$, and the vertices of depth $d$ are precisely the leaves. Thus all the leaves lie at depth $d$. The \emph{height} of a vertex is the number of edges from it to any leaf below it, so that a vertex of depth $j$ has height $d-j$; in particular the root has height $d$ and the leaves have height zero. The \emph{subtree below} a vertex consists of that vertex and all its descendants.

Write the children of a vertex $P$ as $P_1,\ldots,P_M$, and denote the edge from $P$ to $P_i$ by $(P,i)$. There are
\[
\kappa=M+M^2+\cdots+M^d
\]
edges.

To order them, give each edge $(P,i)$ an \emph{address}, namely the sequence of child choices from the root to its endpoint $P_i$ farther from the root. An address is therefore a word
\[
(i_1,\ldots,i_k),\qquad i_j\in\{1,\ldots,M\},\quad1\le k\le d,
\]
and each such word is the address of exactly one edge. For addresses $u\ne v$, we write $u<_{\mathrm{lex}}v$ when either $u$ is a proper prefix of $v$, or the entry of $u$ is smaller at the first position where $u$ and $v$ differ. This is lexicographic order with proper prefixes first, and we order the edges by their addresses. The definition applies to arbitrary $M\ge2$ and $d\ge1$. Purely as an illustration, in the small case $M=2$ and $d=2$ (which plays no special role in what follows), the order is
\[
(1),\ (1,1),\ (1,2),\ (2),\ (2,1),\ (2,2).
\]
Equivalently, starting at the root, one lists $(P,1)$, then every edge below $P_1$ in the same manner, then $(P,2)$ and every edge below $P_2$, and so on through child $M$; this is usually called preorder.

Two consequences of this order will be used repeatedly. First, an edge and all the edges below its child are consecutive in this order, since their addresses are exactly the words beginning with the address of that edge; this is what makes the span estimate of lemma~\ref{lem:span} possible. Second, if $u<_{\mathrm{lex}}v$, then the edge with address $u$ is either an ancestor edge of the edge with address $v$, or leads to an earlier child at some vertex on the path to $v$, or lies below such an earlier child. In particular, the ancestor edges of an edge, and the edges leading to earlier children at the vertices on its path, all precede it; this is why stability, in lemma~\ref{lem:stable}, preserves the routing decisions needed in lemma~\ref{lem:route}.

\subsection{Assigning sequence indices to the edges}\label{subsec:indices}

An edge $e$ receives an interval of consecutive integer indices
\[
W_e=\{a_e,a_e+1,\ldots,b_e\},
\]
corresponding to the terms $q^{a_e},q^{a_e+1},\ldots,q^{b_e}$. By the \emph{length} of $W_e$ we mean the number $|W_e|=b_e-a_e+1$ of indices it contains, not the distance between the corresponding points on the real line.

The intervals are placed in the order of the edges described above, with exactly $g$ unused indices between consecutive intervals: if one interval ends at $b$, the next begins at $b+g+1$. Here ``unused'' only means that these indices are not assigned to any edge; the corresponding terms are still terms of the geometric sequence. By the first consequence of the order recorded above, the interval of an edge $e=(P,i)$ is followed immediately, apart from one gap, by the intervals of all the edges below $P_i$, and only then by the interval of the next sibling edge $(P,i+1)$. We call the \emph{group} of $e$ the set of consecutive indices formed by $W_e$ and, when $P_i$ is not a leaf, the following gap together with the intervals and gaps of all the edges below $P_i$. The lengths will be chosen so that the group of $e$ occupies at most twice the length of $W_e$. We first determine the lengths, and only afterwards the starting indices.

\subsubsection*{The lengths}

Every edge leaving a vertex of height $h$ receives an interval of the same length $r_h$, so that
\[
|W_e|=r_h\qquad\text{when $e$ leaves a vertex of height $h$}.
\]
Here $1\le h\le d$, since leaves have no outgoing edges; the given integer $r_0$ is not the length of any interval, but a lower bound imposed on all of them.

We also record how much room a whole subtree takes in the list of indices. For a vertex of height $h$, let $\sigma_h$ be the number of indices from the first index of the first interval assigned to an edge in its subtree to the last index of the last such interval, gaps included. This span does not include the edge entering the vertex, nor any gaps before or after the subtree. All subtrees of the same height have the same shape and use the same lengths, so $\sigma_h$ depends only on $h$. Thus $r_h$ is the length of one outgoing interval, while $\sigma_h$ is the span of an entire subtree.

These quantities are determined from the leaves toward the root, in the order
\[
r_1,\ \sigma_1,\ r_2,\ \sigma_2,\ \ldots,\ r_d,\ \sigma_d.
\]
At height one the children are leaves, so there are no edges below them. We set
\[
r_1=r_0,
\]
and the subtree consists of $M$ intervals separated by $M-1$ gaps, so that
\[
\sigma_1=Mr_1+(M-1)g.
\]

Now let $h\ge2$, suppose $\sigma_{h-1}$ is known, and let $e=(P,i)$ leave a vertex $P$ of height $h$. The interval $W_e$, the following gap, and the subtree of the child $P_i$, which has height $h-1$, occupy a consecutive span of
\[
\underbrace{r_h}_{W_e}+\underbrace{g}_{\text{gap}}+\underbrace{\sigma_{h-1}}_{\text{child subtree}}
\]
indices. For this to be at most $2r_h$ it is enough that $r_h\ge g+\sigma_{h-1}$, that is, the interval on $e$ should be at least as long as everything that follows it within its group. Keeping also the lower bound $r_0$, we define
\begin{equation}\label{eq:length-choice}
r_h=\max\{r_0,g+\sigma_{h-1}\}\qquad(h\ge2),
\end{equation}
which guarantees
\[
r_h+g+\sigma_{h-1}\le2r_h.
\]
The subtree of $P$ consists of $M$ such groups, one for each child, with one further gap between consecutive groups, so that
\[
\sigma_h=M(r_h+g+\sigma_{h-1})+(M-1)g.
\]
All these quantities are finite.

\subsubsection*{The starting indices}

Choose $n_0\ge1$ with $2q^{n_0}\le\frac{1}{4}$, start the first interval at $n_0$, and place the remaining intervals in order with exactly $g$ unused indices between consecutive intervals. Let
\begin{equation}\label{eq:N-def}
{\mathcal N}=\bigcup_eW_e,
\end{equation}
the union being over all the edges of the tree, be the finite set of indices assigned to edges; this is the set of sequence indices used for the local tests of section~\ref{sec:routing}. Since every $n\in{\mathcal N}$ satisfies $n\ge n_0$, we have
\[
0<tq^n\le\frac{1}{4}\qquad\text{for all }t\in[1,2],
\]
so that all the points $x+tq^n$ with $n\in{\mathcal N}$ lie in $(x,x+\frac{1}{4}]$. This will be used when working with periodic grids.

\subsubsection*{An example}

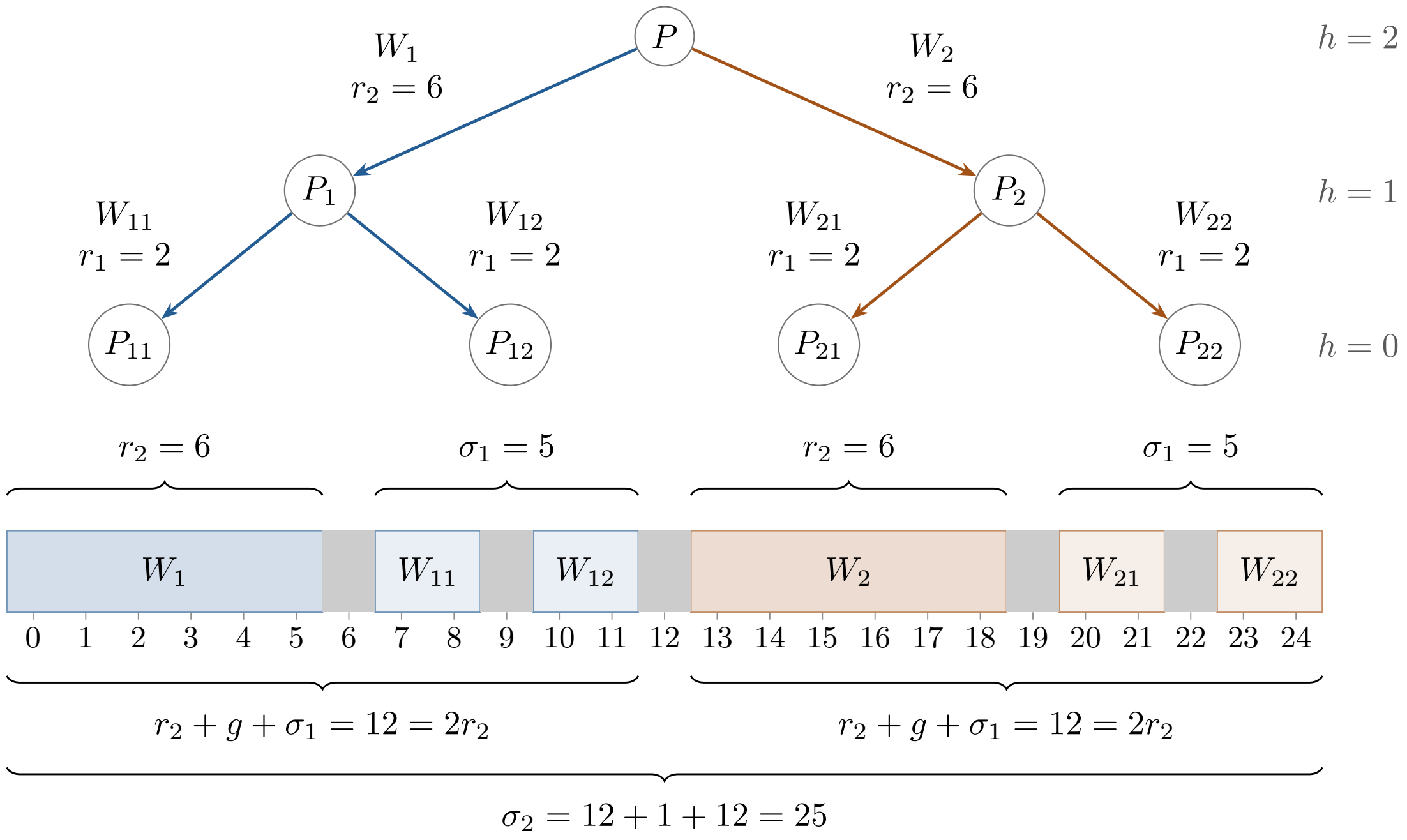
\begin{figure}[htbp]
\centering
\begin{tikzpicture}[
  x=1cm,y=1cm,font=\small,scale=.94,
  vertex/.style={circle,draw=black!55,fill=white,inner sep=2.5pt},
  edge/.style={-{Stealth[length=1.8mm]},line width=.9pt},
  brace/.style={decorate,decoration={brace,amplitude=4pt},line width=.55pt}
]
\node[vertex] (P) at (7.25,4.95) {$P$};
\node[vertex] (P1) at (3.45,3.25) {$P_1$};
\node[vertex] (P2) at (11.05,3.25) {$P_2$};
\node[vertex] (P11) at (1.35,1.55) {$P_{11}$};
\node[vertex] (P12) at (5.55,1.55) {$P_{12}$};
\node[vertex] (P21) at (8.95,1.55) {$P_{21}$};
\node[vertex] (P22) at (13.15,1.55) {$P_{22}$};
\draw[edge,leftbranch] (P)--(P1);
\draw[edge,rightbranch] (P)--(P2);
\draw[edge,leftbranch] (P1)--(P11);
\draw[edge,leftbranch] (P1)--(P12);
\draw[edge,rightbranch] (P2)--(P21);
\draw[edge,rightbranch] (P2)--(P22);
\node[align=center] at (4.3,4.6) {$W_1$\\$r_2=6$};
\node[align=center] at (10.2,4.6) {$W_2$\\$r_2=6$};
\node[align=center] at (1.3,2.75) {$W_{11}$\\$r_1=2$};
\node[align=center] at (5.6,2.75) {$W_{12}$\\$r_1=2$};
\node[align=center] at (8.9,2.75) {$W_{21}$\\$r_1=2$};
\node[align=center] at (13.2,2.75) {$W_{22}$\\$r_1=2$};
\node[anchor=west,text=black!65] at (14.3,4.95) {$h=2$};
\node[anchor=west,text=black!65] at (14.3,3.25) {$h=1$};
\node[anchor=west,text=black!65] at (14.3,1.55) {$h=0$};
\def\cell{.58}
\foreach \a/\b/\name/\col/\shade in {
  0/6/1/leftbranch/20, 7/9/11/leftbranch/10, 10/12/12/leftbranch/10,
  13/19/2/rightbranch/20, 20/22/21/rightbranch/10, 23/25/22/rightbranch/10}{
  \path[fill=\col!\shade,draw=\col!65,line width=.45pt]
    ({\a*\cell},-1.4) rectangle ({\b*\cell},-.5);
  \node at ({(\a+\b)*\cell/2},-.95) {$W_{\name}$};
}
\foreach \a in {6,9,12,19,22}{
  \path[fill=black!20] ({\a*\cell},-1.4) rectangle ({(\a+1)*\cell},-.5);
}
\foreach \k in {0,...,24}{
  \draw[black!45,line width=.35pt] ({(\k+.5)*\cell},-1.4)--({(\k+.5)*\cell},-1.48);
  \node[font=\footnotesize] at ({(\k+.5)*\cell},-1.68) {$\k$};
}
\draw[brace] (0,-.11)--({6*\cell},-.11) node[midway,above=6pt] {$r_2=6$};
\draw[brace] ({7*\cell},-.11)--({12*\cell},-.11) node[midway,above=6pt] {$\sigma_1=5$};
\draw[brace] ({13*\cell},-.11)--({19*\cell},-.11) node[midway,above=6pt] {$r_2=6$};
\draw[brace] ({20*\cell},-.11)--({25*\cell},-.11) node[midway,above=6pt] {$\sigma_1=5$};
\draw[brace] ({12*\cell},-2.10)--(0,-2.10)
  node[midway,below=7pt] {$r_2+g+\sigma_1=12=2r_2$};
\draw[brace] ({25*\cell},-2.10)--({13*\cell},-2.10)
  node[midway,below=7pt] {$r_2+g+\sigma_1=12=2r_2$};
\draw[brace] ({25*\cell},-3.11)--(0,-3.11)
  node[midway,below=7pt] {$\sigma_2=12+1+12=25$};
\end{tikzpicture}
\caption{The illustrative case $M=2$, $d=2$, $g=1$, $r_0=2$, for which $r_1=2$ and $r_2=6$. Above, the tree, with $h$ the height; $W_{11}$ denotes the interval of the edge with address $(1,1)$, and similarly for the others. Below, the same intervals in the order of their addresses. A label $k$ on the strip denotes the index $n=n_0+k$, not the value $q^n$, and the gray positions are the unused indices, one in each gap. Each group, consisting of an edge interval, one gap and the subtree below the child, occupies at most $2r_2$ indices.}
\label{fig:tree-example}
\end{figure}

The following small case, drawn in figure~\ref{fig:tree-example}, is purely an illustration, the actual parameters being chosen in section~\ref{sec:finish}. Take $M=2$, $d=2$, $g=1$ and $r_0=2$. At height one, the two intervals below each child of the root have length $2$ and are separated by one unused index, so that $r_1=2$ and $\sigma_1=2+1+2=5$; note that this span of the subtree rooted at $P_1$ does not include the interval on the edge entering $P_1$. At height two, $r_2=\max\{2,1+5\}=6$, and the interval on the edge entering $P_1$, the following gap, and the subtree rooted at $P_1$ together occupy $6+1+5=12=2r_2$ indices. The same holds for $P_2$, and the two groups are separated by one unused index, so that $\sigma_2=12+1+12=25$.

\subsubsection*{The span estimate}

\begin{figure}[htbp]
\centering
\begin{tikzpicture}[x=1cm,y=1cm,>=Stealth,font=\small]
\draw[fill=blue!9] (0,0) rectangle (5.1,1.35);
\draw[fill=gray!15] (5.1,0) rectangle (5.7,1.35);
\draw[fill=green!9!white] (5.7,0) rectangle (10.1,1.35);
\node[align=center] at (2.55,.675) {$W_e$\\$r_h$ indices};
\node at (5.4,.675) {$g$};
\node[align=center,text width=4cm] at (7.9,.675) {intervals and gaps below the child\\span $\sigma_{h-1}$};
\node[above] at (0,1.4) {$a_e$};
\node[above] at (5.1,1.4) {$b_e$};
\node[above] at (10.1,1.4) {$b_e^*$};
\draw[<->] (0,-.45)--node[below] {total span at most $2r_h$}(10.1,-.45);
\draw[->] (10.25,.52)--(11.3,.52);
\node[align=center,text width=2.3cm] at (11.2,-.45) {groups of later siblings};
\end{tikzpicture}
\caption{The group of an edge $e$ leaving a vertex of height $h$: the interval $W_e$ of $r_h$ indices, one gap of $g$ indices, and the span $\sigma_{h-1}$ of the subtree below the child. The choice $r_h\ge g+\sigma_{h-1}$ gives $r_h+g+\sigma_{h-1}\le2r_h$, which is lemma~\ref{lem:span}. As a preview of section~\ref{sec:scales}: the tests attached to $e$ consult no grid finer than the one indexed by $b_e^*$, so grids in groups of later siblings, however fine, will not enter the count of cell configurations for these tests.}
\label{fig:span}
\end{figure}
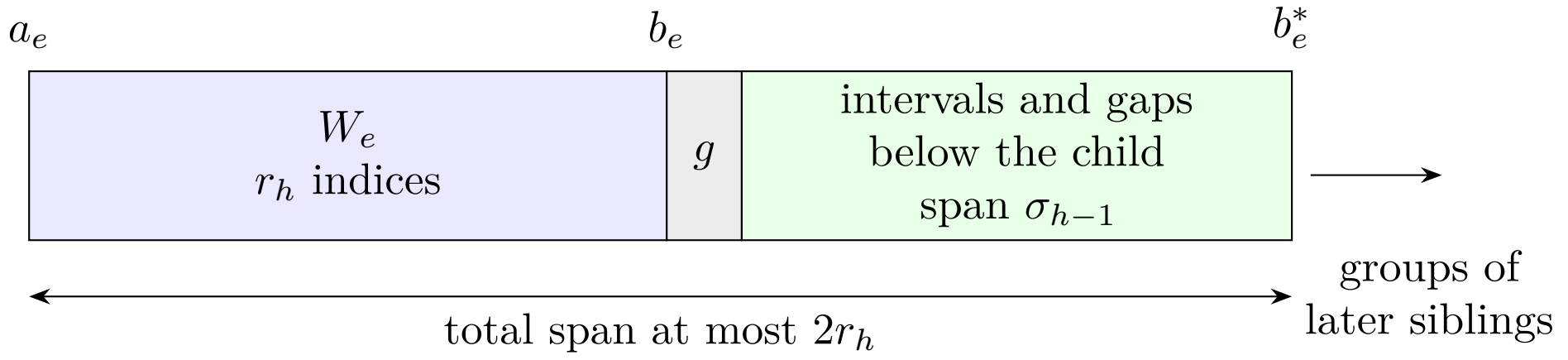

For an edge $e=(P,i)$, which goes from $P$ to its child $P_i$, we say that an edge $f$ lies \emph{below} $P_i$ if it belongs to the subtree below $P_i$, that is, if the endpoint of $f$ nearer the root is $P_i$ or a descendant of $P_i$. These are the edges under $P_i$ at every depth, not only the edges leaving $P_i$. Define
\begin{equation}\label{eq:bstar}
b_e^*=\max\{b_f:\ f=e\ \text{or}\ f\ \text{lies below}\ P_i\},
\end{equation}
the largest index assigned to $e$ or to an edge below $P_i$; if $P_i$ is a leaf, no edge lies below it, and $b_e^*=b_e$. Since the group of $e$ is consecutive, its full span, including the unused gaps, is $\{a_e,\ldots,b_e^*\}$; see figure~\ref{fig:span}. This span has $b_e^*-a_e+1$ indices, of which fewer are actually assigned to edges.

\begin{lemma}\label{lem:span}
If $e$ leaves a vertex of height $h$, then
\begin{equation}\label{eq:span}
b_e^*-a_e+1\le2r_h.
\end{equation}
\end{lemma}

\begin{proof}
If $h=1$, the child is a leaf, the group consists of $W_e$ alone, and its span is $r_1\le2r_1$. If $h\ge2$, the group consists of $W_e$, one gap and the child subtree, so that
\[
b_e^*-a_e+1=r_h+g+\sigma_{h-1}\le2r_h
\]
by~\eqref{eq:length-choice}.
\end{proof}

\subsection{Nested grids for an arbitrary ratio}\label{subsec:grids}

The random set of section~\ref{sec:routing} is built from cells. We divide each unit interval into $N$ equal half-open cells, numbered $0,\ldots,N-1$, and repeat this division periodically on the whole line. Section~\ref{sec:routing} attaches a random entry to each cell number, so that two points lying in cells with the same number consult the same random entry. Selecting periodic cells at random in this way is the random-cell method of Kolountzakis~\cite{Kolountzakis1997}, and nested dyadic intervals of the kind used below also play a central role in Tom's investigation of the dyadic case with Pramanik~\cite{Tom2015}.

Our grids must satisfy two requirements. Their cells must be small enough to separate the points we shall consider (section~\ref{subsec:separation}), and each later grid must refine every earlier one, so that the cell of a point in a fine grid determines its cell in every coarser grid. To satisfy both, set
\[
c=\frac{4}{1-q},
\]
and for each integer $b\ge1$, let $N_b$ be the smallest power of two which is at least $cq^{-b}$, that is,
\[
N_b=2^{\lceil\log_2(cq^{-b})\rceil},
\]
where the ceiling $\lceil u\rceil$ is the smallest integer at least $u$. Thus
\begin{equation}\label{eq:grid-size}
cq^{-b}\le N_b<2cq^{-b}.
\end{equation}
The lower bound controls the cell width, $\frac{1}{N_b}\le\frac{q^b}{c}$, while the choice of powers of two guarantees nesting.

Each grid is described by three objects:
\begin{gridtable}
\begin{tabular}{|l|>{\raggedright\arraybackslash}p{8.6cm}|}
\hline
\emph{object} & \emph{meaning}\\ \hline
$N_b$ & the number of cells in one period\\ \hline
$\Gamma_b=\frac{1}{N_b}{\mathbb Z}$ & the set of cell boundaries\\ \hline
$J_b(z)=\lfloor N_b\{z\}\rfloor$ & the cell label of $z$, that is, the number of its cell modulo one\\ \hline
\end{tabular}
\end{gridtable}
Here $\{z\}=z-\lfloor z\rfloor$ is the fractional part and $\lfloor u\rfloor$ is the greatest integer at most $u$. The cell label belongs to $\{0,\ldots,N_b-1\}$, and boundary points are assigned to the cell on their right. We say ``cell label'' for $J_b(z)$, reserving the word ``address'' for the words describing edges of the tree.

All the $N_b$ are powers of two and are nondecreasing in $b$. Hence if $b'\ge b$, the quotient $\frac{N_{b'}}{N_b}$ is an integer, every cell of the coarser grid is a union of $\frac{N_{b'}}{N_b}$ consecutive cells of the finer grid, and
\begin{equation}\label{eq:nesting}
J_b(z)=\left\lfloor\frac{N_bJ_{b'}(z)}{N_{b'}}\right\rfloor.
\end{equation}
In words, a fine cell label determines every coarser cell label. We use dyadic partitions for this exact nesting even when $q$ is not a dyadic number. Edge $e$ will use the grid indexed by $b_e$.

\begin{figure}[htbp]
\centering
\begin{tikzpicture}[x=1cm,y=1cm,>=Stealth,font=\small]
\foreach \k/\shade in {4/5,5/22,6/5}{
  \path[fill=leftbranch!\shade,draw=black!60] ({(\k-4)*4},1.1) rectangle ({(\k-3)*4},1.8);
  \node at ({(\k-4)*4+2},1.45) {$\k$};
}
\node[anchor=east] at (-.15,1.45) {$N_1=16$};
\foreach \m in {16,...,27}{
  \pgfmathtruncatemacro{\inside}{(\m>=20)&&(\m<=23)}
  \ifnum\inside=1 \path[fill=leftbranch!10,draw=black!60] ({\m-16},0) rectangle ({\m-15},.7);
  \else \path[draw=black!60] ({\m-16},0) rectangle ({\m-15},.7); \fi
  \node[font=\footnotesize] at ({\m-15.5},.35) {$\m$};
}
\path[fill=rightbranch!25,draw=black!60] (7,0) rectangle (8,.7);
\node[font=\footnotesize] at (7.5,.35) {$23$};
\node[anchor=east] at (-.15,.35) {$N_3=64$};
\foreach \a in {4,8}{\draw[densely dashed,black!50] (\a,1.1)--(\a,.7);}
\draw[rightbranch,thick] (7.68,-.25)--(7.68,2.05);
\node[above,text=rightbranch] at (7.68,2.05) {$z=0.37$};
\draw[decorate,decoration={brace,amplitude=4pt,mirror}] (4,-.1)--(8,-.1)
  node[midway,below=5pt] {coarse cell $5$ = fine cells $20,21,22,23$};
\end{tikzpicture}
\caption{Nesting for $q=\frac{1}{2}$, shown near $z=0.37$. Each coarse cell of the grid with $N_1=16$ cells per period consists of four consecutive fine cells of the grid with $N_3=64$. The fine label $J_3(z)=23$ determines the coarse label $J_1(z)=\lfloor\frac{23}{4}\rfloor=5$.}
\label{fig:nested}
\end{figure}
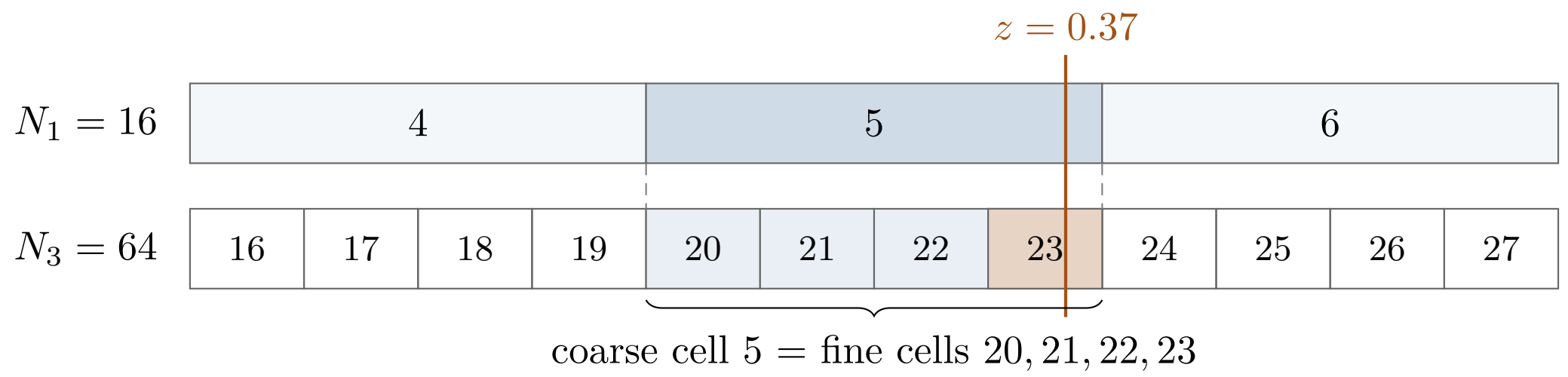

For example, if $q=\frac{1}{2}$, then $c=8$, $N_1=16$ and $N_3=64$. For $z=0.37$,
\[
J_1(z)=\lfloor5.92\rfloor=5,\qquad J_3(z)=\lfloor23.68\rfloor=23,
\]
and indeed $\lfloor\frac{23}{4}\rfloor=5$: the fine cell $23$ is one of the four fine cells $20,21,22,23$ making up the coarse cell $5$, as in figure~\ref{fig:nested}. The point $1.37$ has the same cell labels, because the construction repeats in every unit interval.

\subsection{Small displacements preserve earlier cell labels}\label{subsec:stable}

Fix an edge $e$ and an index $n\in W_e$. We ask when the base point $x$ and the point $x+tq^n$ lie in cells with the same label in every grid used by an edge preceding $e$. Three observations reduce this question to a single condition.
\begin{itemize}
\item Since $n\ge a_e$ and $t\le2$, the displacement satisfies $0<tq^n\le2q^{a_e}$, and it is to the right.
\item The intervals are placed in increasing order, so if $b'$ denotes the last index of the interval immediately before $W_e$, then $b_f\le b'$ for every edge $f$ preceding $e$. The grid indexed by $b'$ therefore refines all the earlier grids, and by~\eqref{eq:nesting} it suffices to check that grid.
\item The two points have the same label in that grid if moving to the right from $x$ by any distance up to $2q^{a_e}$ does not reach the next boundary of $\Gamma_{b'}$.
\end{itemize}
Recall that $a_e=b'+g+1$. We call a base point $x$ \emph{stable} if, for every interval $W_e$ other than the first,
\begin{equation}\label{eq:stable-definition}
(x,x+2q^{a_e}]\cap\Gamma_{b'}=\varnothing.
\end{equation}
Let ${\mathcal S}$ be the set of stable base points. A boundary at $x$ itself is allowed, since boundary points belong to the cell on their right and the movement starts inside that cell; reaching a boundary is not allowed, since the moving point would then enter a new cell.

The role of the gaps is now visible. The earlier cells have width $1/N_{b'}$, so for every $n\in W_e$ and $t\in[1,2]$, \eqref{eq:grid-size} gives
\[
\frac{tq^n}{1/N_{b'}}=tq^nN_{b'}\le2q^{a_e}N_{b'}\le4cq^{a_e-b'}=4cq^{g+1}.
\]
The left side is the ratio of the displacement to the width of the earlier cells.
Increasing $g$ makes the possible displacement small compared with the earlier cells, and this is the geometric purpose of leaving unused indices between the intervals. The base points excluded by~\eqref{eq:stable-definition} are easy to describe: immediately to the left of each boundary of $\Gamma_{b'}$ lies an interval of length $2q^{a_e}$ of starting points which could reach that boundary. The measure estimate below counts these boundaries and multiplies by this length.

\begin{lemma}\label{lem:stable}
The set ${\mathcal S}$ is measurable and $1$-periodic, with
\begin{equation}\label{eq:unstable-size}
\rho({\mathbb R}\setminus{\mathcal S})\le4\kappa cq^{g+1},
\end{equation}
where $\kappa=M+M^2+\cdots+M^d$ is the number of edges, which is also the number of index intervals. If $x\in{\mathcal S}$, $n\in W_e$, and $t\in[1,2]$, then every edge $f$ preceding $e$ satisfies
\begin{equation}\label{eq:earlier}
J_{b_f}(x+tq^n)=J_{b_f}(x).
\end{equation}
\end{lemma}

\begin{proof}
Fix an interval other than the first, and write $\Delta=2q^{a_e}$. A base point $x$ violates~\eqref{eq:stable-definition} exactly when some boundary $\gamma\in\Gamma_{b'}$ satisfies $\gamma\in(x,x+\Delta]$, and
\[
\gamma\in(x,x+\Delta]\iff x<\gamma\le x+\Delta\iff\gamma-\Delta\le x<\gamma.
\]
Hence the base points violating~\eqref{eq:stable-definition} are precisely the points of
\[
\bigcup_{\gamma\in\Gamma_{b'}}[\gamma-\Delta,\gamma).
\]
Only intervals to the left of the boundaries occur: a base point just to the right of $\gamma$ moves away from $\gamma$ as the displacement grows, so it cannot reach that boundary. For example, if $\gamma=1$ and $\Delta=0.1$, the base points which reach $\gamma$ form $[0.9,1)$. There are $N_{b'}$ boundaries per period, so these intervals have total density at most
\[
N_{b'}\Delta=2N_{b'}q^{a_e}\le4cq^{a_e-b'}=4cq^{g+1}.
\]
Each of the $\kappa$ index intervals other than the first gives one such condition, so there are $\kappa-1$ conditions. Adding over them proves~\eqref{eq:unstable-size}, and also shows measurability and periodicity.

For a stable $x$, the displacement satisfies $0<tq^n\le2q^{a_e}$, so no boundary of $\Gamma_{b'}$ is crossed or reached and $J_{b'}(x+tq^n)=J_{b'}(x)$. Since $b_f\le b'$ for every earlier edge, the nesting formula~\eqref{eq:nesting} proves~\eqref{eq:earlier}. For the first interval there is nothing to prove.
\end{proof}

\subsection{Distinct cell labels at the current resolution}\label{subsec:separation}

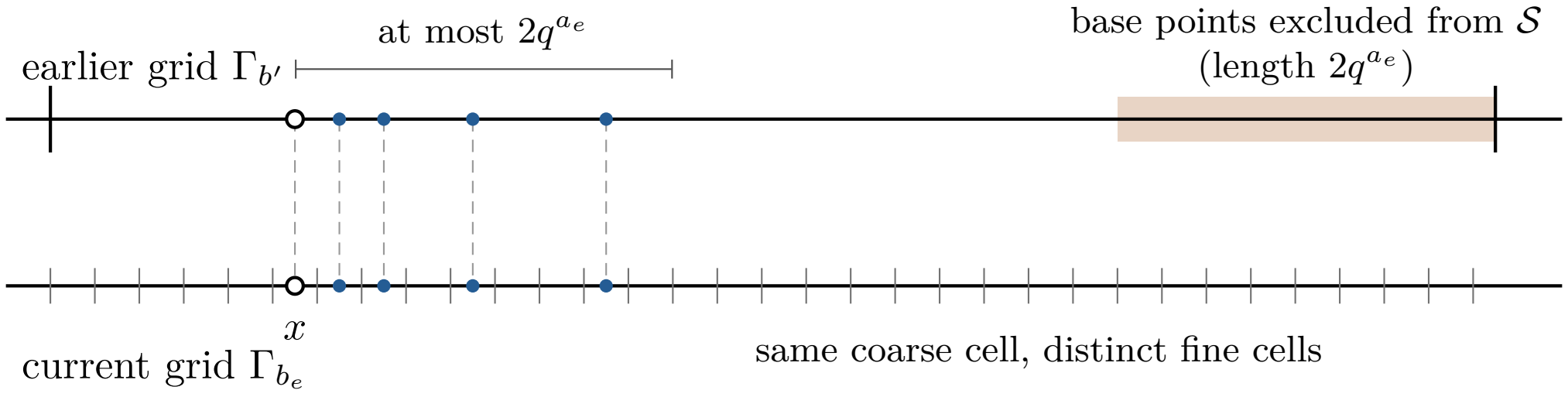
\begin{figure}[htbp]
\centering
\begin{tikzpicture}[x=1cm,y=1cm,>=Stealth,font=\small,scale=.95]
\node[anchor=west] at (-.4,1.95) {earlier grid $\Gamma_{b'}$};
\fill[rightbranch!25] (9.6,1.3) rectangle (13,1.7);
\draw[thick] (-.4,1.5)--(13.6,1.5);
\foreach \a in {0,13}{\draw[thick] (\a,1.2)--(\a,1.8);}
\node[align=center,font=\footnotesize] at (11.3,2.15) {base points excluded from ${\mathcal S}$\\(length $2q^{a_e}$)};
\node[anchor=west] at (-.4,-.75) {current grid $\Gamma_{b_e}$};
\draw[thick] (-.4,0)--(13.6,0);
\foreach \k in {0,...,32}{\draw[black!55] ({\k*.4},-.16)--({\k*.4},.16);}
\foreach \a in {2.6,3.0,3.8,5.0}{
  \draw[densely dashed,black!40] (\a,1.5)--(\a,0);
  \fill[leftbranch] (\a,1.5) circle (.06);
  \fill[leftbranch] (\a,0) circle (.06);
}
\draw[densely dashed,black!40] (2.2,1.5)--(2.2,0);
\draw[black,fill=white,thick] (2.2,1.5) circle (.075); \draw[black,fill=white,thick] (2.2,0) circle (.075);
\node[below] at (2.2,-.18) {$x$};
\draw[|-|,black!70] (2.2,1.95)--(5.6,1.95);
\node[above,font=\footnotesize] at (3.9,1.98) {at most $2q^{a_e}$};
\node[align=left,font=\footnotesize,anchor=west] at (6.2,-.6) {same coarse cell, distinct fine cells};
\end{tikzpicture}
\caption{The same points $x$ (open circle) and $x+tq^n$, $n\in W_e$ (filled), drawn on an earlier grid and on the grid of the current edge, with vertical guides joining them. The displacements are at most $2q^{a_e}$, so a stable $x$ and all these points share one coarse cell, and the earlier cell labels, hence the earlier routing decisions, are unchanged. In the current grid the same points occupy distinct cells and will consult distinct random entries. A base point in the shaded interval to the left of a coarse boundary could reach that boundary and is therefore not stable; a base point lying exactly on a boundary is allowed. The picture is schematic.}
\label{fig:grids}
\end{figure}

The previous lemma asks nearby points to share their labels in the earlier grids, whereas in the grid of the current edge we need the opposite property. For a stable base point $x$ and an edge $e$, the two facts hold simultaneously:
\begin{gridtable}
\begin{tabular}{|>{\raggedright\arraybackslash}p{4.0cm}|>{\raggedright\arraybackslash}p{8.6cm}|}
\hline
\emph{grid} & \emph{what happens to $x$ and the points $x+tq^n$, $n\in W_e$}\\ \hline
an earlier grid & they all have the same cell label\\ \hline
the grid indexed by $b_e$ & they have pairwise distinct cell labels\\ \hline
\end{tabular}
\end{gridtable}
They are compatible because one coarse cell contains many fine cells; see figure~\ref{fig:grids}. Since $t>0$ and the powers of $q$ decrease, the points are ordered as
\[
x<x+tq^{b_e}<x+tq^{b_e-1}<\cdots<x+tq^{a_e}.
\]
Increasing the exponent moves the point toward $x$; we state this because the indices themselves are listed in increasing order.

Cell labels are defined modulo one, so we must measure distances modulo one as well. We use the circle ${\mathbb R}/{\mathbb Z}$, on which the distance between the images of $u$ and $v$ is $\min_{k\in{\mathbb Z}}|u-v-k|$. For example, $0.99$ and $1.01$ are close, although their fractional parts lie near opposite ends of $[0,1)$, and the circle records this correctly. Since all our points lie in an interval of length at most $\frac{1}{4}$, their circular distances agree with their ordinary distances.

\begin{lemma}\label{lem:separation}
For every edge $e$, base point $x$, and $t\in[1,2]$, the points
\[
\{x\}\cup\{x+tq^n:n\in W_e\}
\]
have pairwise distinct cell labels under $J_{b_e}$. They also have pairwise distinct cell labels in every finer grid.
\end{lemma}

\begin{proof}
Write $b=b_e$. The proof compares two lengths. On the one hand, the minimum separation of the points is at least $(1-q)q^b$. Indeed, the distance from $x$ to the nearest of the other points is $tq^b\ge q^b$, and if $n<n'$ belong to $W_e$, then
\[
t(q^n-q^{n'})=tq^n(1-q^{n'-n})\ge(1-q)q^b.
\]
On the other hand, by~\eqref{eq:grid-size} the cell width satisfies
\[
\frac{1}{N_b}\le\frac{q^b}{c}=\frac{(1-q)q^b}{4}.
\]
All the points lie in an interval of length at most $2q^{n_0}\le\frac{1}{4}$, so their circular distances equal their ordinary distances. Indeed, if $|u-v|\le\frac{1}{2}$, then $|u-v-k|\ge|k|-|u-v|\ge\frac{1}{2}\ge|u-v|$ for every integer $k\ne0$, so the minimum defining the circular distance is attained at $k=0$. 

Two points with the same cell label have circular distance strictly less than the cell width, which the separation rules out. Distinctness persists in finer grids by~\eqref{eq:nesting}, since equal fine labels would force equal coarse labels.
\end{proof}

\section{Random tables and the set they define}\label{sec:routing}

As in Kolountzakis's random-cell method~\cite{Kolountzakis1997}, the small set $B$ will be a union of cells chosen at random. In the OpenAI construction, a point first travels down the tree, and only the table at the leaf it reaches decides whether the point belongs to $B$. We describe the construction in four stages.

\subsection{Two kinds of random tables}

For each edge $e=(P,i)$ with $i<M$, attach a random routing table
\[
S_e:\{0,\ldots,N_{b_e}-1\}\longrightarrow\{0,1\},
\]
that is, a family $\bigl(S_e(j)\bigr)_{j=0}^{N_{b_e}-1}$ of independent Bernoulli-$\frac{1}{2}$ random variables, one for each cell label. We call $S_e(j)$ the \emph{entry} of $S_e$ at the label $j$, and we call the entries of the routing tables \emph{routing entries}.
When this table is consulted for a point $z$, the value $S_e(J_{b_e}(z))=1$ sends the point along $e$ to the child $P_i$, whereas the value $0$ means that the next child $P_{i+1}$ is considered. There is no routing table on the last edge $(P,M)$, since the last child is chosen by default when the others have all been rejected. For each leaf $L$, let $b(L)$ be the last index of the interval on the edge entering $L$, and attach a terminal table
\[
T_L:\{0,\ldots,N_{b(L)}-1\}\longrightarrow\{0,1\}
\]
whose entries $T_L(j)$, called \emph{terminal entries}, are independent Bernoulli-$p$ variables.
All the entries in all the tables are mutually independent, and since the tree and every table are finite, this is a finite probability space. Below, to \emph{reveal} some entries means to condition on their values: we hold those values fixed and compute probabilities by averaging over the remaining entries. This terminology is standard in probabilistic arguments of this kind. We use three further terms.
\begin{itemize}
\item An entry is identified by its table and its cell label; for instance, the pair $(L,j)$ identifies the terminal entry $T_L(j)$, which is a random variable taking the values zero and one. When some information determines which entry is read, it determines the pair $(L,j)$, not the value of $T_L(j)$.
\item A \emph{realization} of the tables is an assignment of a value, zero or one, to every entry of every table. There are finitely many realizations.
\item Once the values of some entries have been fixed, an assignment of values to all the remaining entries is called a \emph{completion}.
\end{itemize}

All the table entries are sampled once and for all. Every point consulting a given entry sees the same value; no new coin is tossed when another point arrives. The two kinds of entries have the following meanings:
\begin{gridtable}
\begin{tabular}{|>{\raggedright\arraybackslash}p{4.2cm}|>{\raggedright\arraybackslash}p{4.0cm}|>{\raggedright\arraybackslash}p{4.0cm}|}
\hline
\emph{table} & \emph{entry $1$} & \emph{entry $0$}\\ \hline
routing table $S_e$, $e=(P,i)$ & choose child $i$ & consider the next child\\ \hline
terminal table $T_L$ & include the point in $B$ & exclude the point from $B$\\ \hline
\end{tabular}
\end{gridtable}
A point $z$ consults the entry of a table whose index is its cell label in that table's grid. The tree describes which entries are consulted; the real number $z$ itself does not change while its route is being determined.

\subsection{Following a point through the tree}

Given $z\in{\mathbb R}$, start at the root. At a vertex $P$, read
\[
S_{(P,1)}(J_{b_{(P,1)}}(z)),\ldots,
S_{(P,M-1)}(J_{b_{(P,M-1)}}(z))
\]
in order, and choose the first child whose entry is one, taking child $M$ if every entry is zero. For instance, when $M=3$:
\begin{gridtable}
\begin{tabular}{|c|c|c|}
\hline
\emph{first routing entry} & \emph{second routing entry} & \emph{child chosen}\\ \hline
$1$ & not consulted & $1$\\ \hline
$0$ & $1$ & $2$\\ \hline
$0$ & $0$ & $3$, by default\\ \hline
\end{tabular}
\end{gridtable}
The choice of a child at a vertex is called a \emph{routing decision}. Taking the last child $M$ because all the other children have been rejected is called a \emph{default choice}. If the route of a point makes at least one default choice, the first vertex at which it does so is called its \emph{first default vertex}; at that vertex all the $M-1$ routing entries read by the point are zero, and the route continues to the last child.

Continue in the same way from the chosen child until a leaf is reached, call this leaf $L(z)$, and define the random set
\begin{equation}\label{eq:random-set}
B=\{z:T_{L(z)}(J_{b(L(z))}(z))=1\}.
\end{equation}
Thus, once the route reaches the leaf $L(z)$, we inspect the value $T_{L(z)}(J_{b(L(z))}(z))$: if it is one, the point $z$ belongs to $B$, and if it is zero, $z$ does not belong to $B$. In either case the procedure stops at this leaf; a zero does not send the point back to an earlier vertex to try another branch. This is the difference between the two kinds of zeros: a routing zero leads to the next child, while a terminal zero decides exclusion from $B$. For a fixed realization of the tables, all the decisions are constant on the cells of the finest grid in the construction, so $B$ is $1$-periodic and is a finite union of half-open cells in each period.

\begin{lemma}\label{lem:density}
The expected density of $B$ is $p$.
\end{lemma}

\begin{proof}
Fix $z$ and reveal first every routing table. This determines the leaf $L(z)$ and which entry of $T_{L(z)}$ is read, namely the one with cell label $J_{b(L(z))}(z)$, while the value of that entry is still one with probability $p$. Since every point reaches exactly one leaf, and the leaf reached depends only on the routing tables, which are independent of the terminal tables,
\[
{\mathbb P}(z\in B)=\sum_{L\text{ leaf}}{\mathbb P}\bigl(L(z)=L,\ T_L(J_{b(L)}(z))=1\bigr)
=p\sum_{L\text{ leaf}}{\mathbb P}\bigl(L(z)=L\bigr)=p.
\]

Let $\omega$ range over the finitely many realizations of the tables, write ${\mathbb P}(\omega)$ for the probability of $\omega$, and write $B_\omega$ for the set~\eqref{eq:random-set} determined by $\omega$. Then
\[
\begin{aligned}
\mathbb E\,\rho(B)
&=\sum_{\omega}\mathbb P(\omega)
\int_0^1\mathbf 1_{B_\omega}(z)\,dz\\
&=\int_0^1\sum_{\omega}\mathbb P(\omega)
\mathbf 1_{B_\omega}(z)\,dz\\
&=\int_0^1\mathbb P(z\in B)\,dz=p.
\end{aligned}
\]
The second equality is only an interchange of an integral with a finite weighted sum.
\end{proof}

The gaps of $g$ unused indices between the index intervals play no role here. They only affect which points of a copy are used in the later local tests, and they do not exclude any real point from the definition of $B$: every real point $z$ follows a route and reads a terminal entry, so ${\mathbb P}(z\in B)=p$ for every $z$, including the points $x+tq^n$ whose indices $n$ lie in the gaps.

Individual realizations need not have density exactly $p$; only the expectation is controlled, which suffices for section~\ref{sec:finish}. This calculation is the reason for keeping the terminal tables separate: a tree with many branches retains expected density $p$, and branching does not require us to select a separate set of density $p$ for every branch and take their union.

\subsection{Fixing the routing entries at the base point}

Fix a base point $x$. Before anything is revealed, its route is random. When the route reaches a vertex, the last child is taken by default exactly when all $M-1$ routing entries there are zero, which has probability $2^{-(M-1)}$. The entries read at a newly reached vertex belong to tables not used at its ancestors, so revealing the route one vertex at a time proves
\begin{equation}\label{eq:no-default}
{\mathbb P}(\text{the route of }x\text{ never makes a default choice})
=(1-2^{-(M-1)})^d.
\end{equation}

Now, in every routing table $S_e$, consider the single entry at the cell label of $x$, that is, $S_e(J_{b_e}(x))$. Specify a value $\xi_e\in\{0,1\}$ for each of these entries, write $\xi=(\xi_e)_e$, and condition on the event
\begin{equation}\label{eq:specification}
S_e(J_{b_e}(x))=\xi_e\qquad\text{for every edge }e=(P,i)\text{ with }i<M.
\end{equation}
We call $\xi$ a \emph{specification} of the routing entries at $x$, or simply a specification, and we abbreviate conditioning on this event as conditioning on $\xi$; every specification has positive probability. The specification fixes exactly one entry in each routing table, including the tables on edges outside the route followed by $x$, and these values determine the entire route of $x$. All the other routing entries and all the terminal entries remain random. Since all the table entries were chosen independently, conditioning on $\xi$ leaves the remaining routing entries independent Bernoulli-$\frac{1}{2}$ random variables, and leaves the terminal entries independent with their original Bernoulli-$p$ distributions. We say that an entry is \emph{unrevealed by} $\xi$ if its value is not fixed by $\xi$, in accordance with the terminology introduced at the beginning of this section, where revealing entries means conditioning on their values. Thus, conditionally on $\xi$, an unrevealed entry keeps its original distribution, and distinct unrevealed entries are independent.

\subsection{Why a default vertex is useful}

\begin{figure}[htbp]
\centering
\begin{tikzpicture}[>=Stealth,font=\small,scale=.95,
  vertex/.style={circle,draw=black!60,fill=white,inner sep=2pt}]
\node[vertex] (R) at (0,4.6) {root};
\node[vertex] (u) at (0,2.6) {$U$};
\node[vertex] (a) at (-5.5,.5) {$U_1$};
\node[vertex] (b) at (0,.5) {$U_2$};
\node[vertex] (c) at (5.5,.5) {$U_3$};
\draw[leftbranch,line width=1.1pt] ($(R.south)+(-.07,0)$)--($(u.north)+(-.07,0)$);
\draw[rightbranch,line width=1.1pt,dash pattern=on 3.5pt off 2pt] ($(R.south)+(.07,0)$)--($(u.north)+(.07,0)$);
\node[anchor=east,font=\footnotesize] at (-.12,4.0) {$x$};
\node[anchor=west,font=\footnotesize] at (.12,4.0) {$y$};
\node[anchor=west,align=left,font=\footnotesize] at (.3,3.5) {same choices for $x$ and $y$\\(same earlier cell labels, by stability)};
\draw[->,black!50] (u)--node[midway,sloped,above,font=\footnotesize,text=black]{$S_{e_1}=0$ at $x$ and at $y$}(a);
\draw[->,rightbranch,line width=1.1pt,dash pattern=on 3.5pt off 2pt] (u)--(b);
\node[anchor=west,align=left,font=\footnotesize,text width=2.6cm] at (.35,.95) {$S_{e_2}=0$ at $x$;\\unrevealed entry\\at $y$, here $1$};
\draw[->,leftbranch,line width=1.1pt] (u)--node[midway,sloped,above,font=\footnotesize]{last child, taken by $x$}(c);
\draw[black!40] (a)--(-6.5,-1.6)--(-4.5,-1.6)--cycle;
\draw[black!40] (b)--(-1,-1.6)--(1,-1.6)--cycle;
\draw[black!40] (c)--(4.5,-1.6)--(6.5,-1.6)--cycle;
\draw[rightbranch,line width=1.1pt,dash pattern=on 3.5pt off 2pt] (b.south)--(.35,-.55)--(.2,-1.6);
\fill[rightbranch] (.2,-1.6) circle (.07);
\draw[leftbranch,line width=1.1pt] (c.south)--(5.2,-.55)--(5.9,-1.6);
\fill[leftbranch] (5.9,-1.6) circle (.07);
\node[align=center,font=\footnotesize,text width=3.4cm] at (-5.5,-2.3) {not entered by $y$};
\node[align=center,font=\footnotesize,text width=5.4cm] at (0,-2.4) {leaf $L_2(y)$: if its terminal entry is $1$, then $Q_2(y)=1$ and $y\in B$};
\node[align=center,font=\footnotesize,text width=3.4cm] at (5.5,-2.3) {route of $x$};
\end{tikzpicture}
\caption{A first default vertex $U$ with $M=3$. The route of the base point $x$ is drawn solid, and that of a point $y=x+tq^n$, $n\in W_{e_2}$, dashed. By stability, $y$ follows the route of $x$ down to $U$ and inherits the rejection of child $1$. At $e_2$ its cell label differs from that of $x$, so it reads an unrevealed routing entry, and if this entry is one it enters $U_2$ while $x$ takes the last child. For points with $n\in W_{e_1}$, the routing entry at $e_1$ itself is unrevealed and no earlier child must be rejected. Points attached to different children use terminal tables in disjoint subtrees.}
\label{fig:default}
\end{figure}
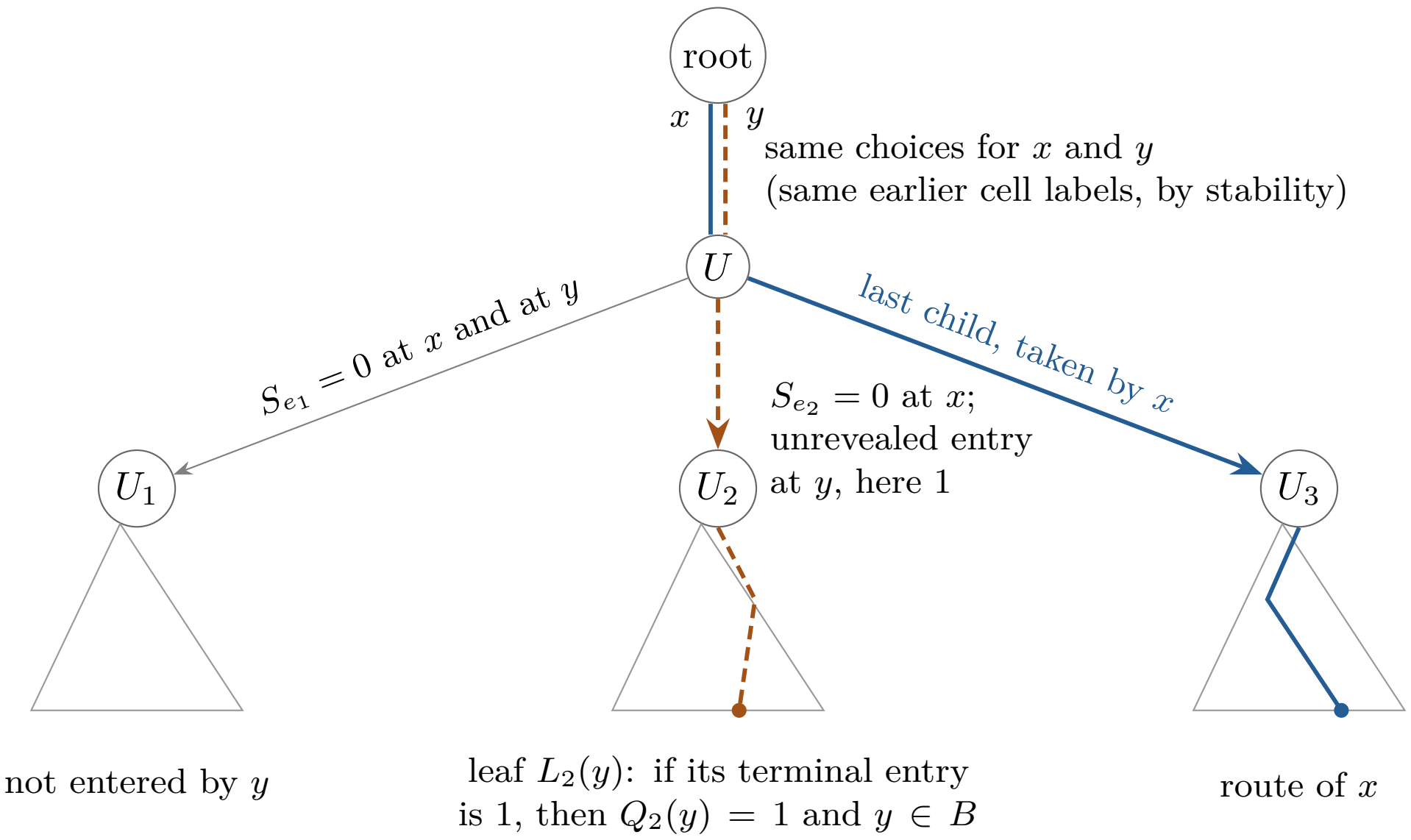

Suppose that $x$ is stable and that the route of $x$ determined by $\xi$ makes at least one default choice, and let $U$ be its first default vertex. Thus all the $M-1$ routing entries read by $x$ at $U$ are zero, and the route of $x$ continues to the last child $U_M$. The other children of $U$ nevertheless remain available to the points of the copy attached to the edges leaving $U$, as the following example with $M=3$ shows; see figure~\ref{fig:default}. Take $y=x+tq^n$ with $n\in W_{(U,2)}$ and $t\in[1,2]$.
\begin{itemize}
\item The edges consulted on the route from the root to $U$ precede $(U,2)$, so by lemma~\ref{lem:stable} the point $y$ reads the same entries as $x$ there and also reaches $U$.
\item The edge $(U,1)$ also precedes $(U,2)$, so $y$ reads the same zero as $x$ there, and child $1$ is rejected for $y$ as well.
\item At $(U,2)$, lemma~\ref{lem:separation} gives $y$ a cell label different from that of $x$, so $y$ consults a routing entry which $\xi$ has not revealed; conditionally on $\xi$, this entry is one with probability $\frac{1}{2}$.
\item If that entry is one, then $y$ enters child $2$, even though $x$ takes the last child $3$.
\end{itemize}
In the same way, for every $i<M$, the points whose indices lie in $W_{(U,i)}$ reach $U$, reject the children $U_1,\ldots,U_{i-1}$ as $x$ does, and then read an unrevealed entry on $(U,i)$, which may send them into $U_i$. This is where the default choice helps. At a vertex $P$ where the route of $x$ chooses a child $P_a$ with $a<M$, a point attached to an edge $(P,i)$ with $i>a$ reads on $(P,a)$ the same entry as $x$, which is one, and follows $x$ into $P_a$, so only the children $P_1,\ldots,P_a$ offer such chances. At a default vertex, $x$ reads zeros on all the edges $(U,1),\ldots,(U,M-1)$, so every non-last child of $U$ offers additional chances, decided by unrevealed entries.

We now set up the general notation. If $U$ has height $h$, put $r=r_h$ and $e_i=(U,i)$ for $1\le i<M$, so that all the routing entries at the cell label of $x$ on these edges are zero. We shall test the points
\begin{equation}\label{eq:tested-points}
y_n(t)=x+tq^n,\qquad n\in W_{e_i},\ 1\le i<M,
\end{equation}
and there are $(M-1)r$ such tests in total; we call these points the \emph{tested points}. For a point $z$, let $L_i(z)$ be the leaf obtained as follows: start at the child $U_i$ and follow the same routing rule, whether or not the actual route of $z$ enters $U_i$. Define the \emph{local test}
\begin{equation}\label{eq:local-test}
Q_i(z)=S_{e_i}(J_{b_{e_i}}(z))
T_{L_i(z)}(J_{b(L_i(z))}(z)).
\end{equation}
We say that this local test succeeds if $Q_i(z)=1$, and fails if $Q_i(z)=0$. Thus success means that the routing entry on $e_i$ accepts $z$ and that the terminal entry reached below $U_i$ is one; only the edge $e_i$ and the subtree below $U_i$ are consulted. Under the assumptions of the next lemma, a successful local test guarantees membership in $B$; a failed local test does not guarantee exclusion from $B$, since the point may still belong to $B$ for another reason.

\begin{lemma}\label{lem:route}
Let $x$ be a stable base point, and let $\xi$ be a specification of the routing entries at $x$ for which the route of $x$ makes at least one default choice. Let $U$ be the first default vertex of this route, put $e_i=(U,i)$ for $1\le i<M$, and let $Q_i$ be the local tests defined in~\eqref{eq:local-test}. Fix a realization of the tables in which every routing table $S_e$ takes the value $\xi_e$ at the cell label of $x$, that is, $S_e(J_{b_e}(x))=\xi_e$. Then, for every $i<M$, every $n\in W_{e_i}$ and every $t\in[1,2]$, the route of the point $y_n(t)=x+tq^n$ reaches $U$ and rejects the children $U_1,\ldots,U_{i-1}$. In particular,
\begin{equation}\label{eq:test-hit}
Q_i(y_n(t))=1\Longrightarrow y_n(t)\in B.
\end{equation}
\end{lemma}

\begin{proof}
We follow the route of $x$ from the root. Let $P$ be a vertex of this route lying strictly above $U$, and suppose that the route of $x$ goes from $P$ to the child $P_a$. Since $U$ is the first default vertex, the choice at $P$ is not a default choice, so $a<M$, the routing entries read by $x$ on the edges $(P,1),\ldots,(P,a-1)$ are zero, and the one read on $(P,a)$ is one. The edge $(P,a)$ is an ancestor edge of $e_i$, and the edges $(P,1),\ldots,(P,a-1)$ lead to earlier children at a vertex on the path to $e_i$; by the second consequence of the edge order recorded in section~\ref{sec:geometry}, all of them precede $e_i$. By lemma~\ref{lem:stable}, the cell labels of $y_n(t)$ in the grids of these edges agree with those of $x$, so at $P$ the point $y_n(t)$ reads the same routing entries as $x$, finds the same values, and also moves to $P_a$. Applying this at each vertex of the route of $x$ above $U$, starting at the root, shows that $y_n(t)$ follows the route of $x$ down to $U$.

At $U$, the edges $e_j$ with $j<i$ lead to earlier children at the vertex where $e_i$ starts, so they also precede $e_i$. By lemma~\ref{lem:stable} again, on each of them $y_n(t)$ reads the same entry as $x$, namely $S_{e_j}(J_{b_{e_j}}(x))=\xi_{e_j}$, which is zero because the route of $x$ makes a default choice at $U$. Hence the children $U_1,\ldots,U_{i-1}$ are rejected.

Finally, suppose that $Q_i(y_n(t))=1$, so that both factors in~\eqref{eq:local-test} equal one. The first factor is the routing entry read by $y_n(t)$ on $e_i$, so the actual route of $y_n(t)$ enters $U_i$. From there it follows the routing rule, exactly as in the definition of $L_i(y_n(t))$, so it ends at the leaf $L(y_n(t))=L_i(y_n(t))$. The second factor is therefore the terminal entry read by $y_n(t)$ at its leaf, and it equals one, so $y_n(t)\in B$ by~\eqref{eq:random-set}.
\end{proof}

\section{The probability that all the local tests fail}\label{sec:independence}

Throughout this section we fix
\begin{itemize}
\item a stable base point $x$;
\item a specification $\xi$ of the routing entries at $x$, for which the route of $x$ has a first default vertex $U$;
\item one dilation $t\in[1,2]$.
\end{itemize}
Consequently the vertex $U$, its height $h$, the length $r=r_h$, the edges $e_i=(U,i)$, and the $(M-1)r$ tested points $y_n(t)$ of~\eqref{eq:tested-points} are all fixed, while every table entry not revealed by $\xi$ is still random. This section concerns one fixed dilation; the passage to all dilations is the subject of section~\ref{sec:scales}.

\subsection{The answer for one tested point}

Consider a single tested point $y=y_n(t)$ with $n\in W_{e_i}$. Its local test $Q_i(y)$ fails either because its routing entry on $e_i$ rejects it, or because that routing entry accepts it and the terminal entry it then reaches is zero. The routing entry is an unrevealed Bernoulli-$\frac{1}{2}$ random variable, and the terminal entry is one with probability $p$, so the failure probability is
\[
\underbrace{\frac{1}{2}}_{\text{routing entry rejects}}
+\underbrace{\frac{1}{2}(1-p)}_{\text{routing entry accepts, terminal entry rejects}}
=1-\frac{p}{2}.
\]
The following lemma states that the failure probabilities of all $(M-1)r$ tested points can be multiplied. This multiplication is the familiar computation at the heart of Kolountzakis's random-cell method~\cite{Kolountzakis1997}, in which $N$ points lying in distinct random cells all miss the selected cells with probability $(1-p)^N$. What requires justification here is that the tested points really read distinct, unrevealed entries.

\begin{lemma}\label{lem:failure}
Fix a stable base point $x$ and a specification $\xi$ for which the first default vertex $U$ exists. Let $r=r_h$, where $h$ is its height, let $e_i=(U,i)$ for $1\le i<M$, and let $y_n(t)$ and $Q_i$ be the tested points~\eqref{eq:tested-points} and the local tests~\eqref{eq:local-test}. For each fixed $t\in[1,2]$,
\begin{equation}\label{eq:exact-failure}
{\mathbb P}(Q_i(y_n(t))=0\text{ for every }i<M,\ n\in W_{e_i}\mid\xi)
=\Bigl(1-\frac{p}{2}\Bigr)^{(M-1)r}.
\end{equation}
\end{lemma}

The rest of this section explains and proves this lemma.

\subsection{Why the argument proceeds in stages}

The routing entries read by the tested points are distinct. Tested points attached to different children of $U$ use different tables, and tested points attached to the same child have distinct cell labels in every grid at least as fine as the grid of $e_i$, by lemma~\ref{lem:separation}; every edge below $U_i$ uses such a grid. The difficulty lies elsewhere: which terminal entry a tested point reads depends on its route below $U_i$, and that route depends on random routing entries. We therefore reveal the routing entries first. This fixes which terminal entries are read, that is, their leaf tables and cell labels, but not the values stored in them; we then check that these entries are distinct. The argument proceeds as follows.
\begin{gridtable}
\begin{tabular}{|>{\raggedright\arraybackslash}p{3.0cm}|>{\raggedright\arraybackslash}p{5.0cm}|>{\raggedright\arraybackslash}p{5.0cm}|}
\hline
\emph{step} & \emph{what is determined or computed} & \emph{what remains to be averaged}\\ \hline
condition on $\xi$ & the route of $x$, the vertex $U$, and the tested points & the routing values not fixed by $\xi$, and the terminal values\\ \hline
condition also on all remaining routing values & each acceptance variable $A_{i,n}$, and the leaf table and cell label of the terminal entry read by each local test & the terminal values\\ \hline
average over the terminal values & the conditional failure probability $(1-p)^Z$ & the routing values not fixed by $\xi$\\ \hline
average over those routing values & the failure probability $(1-\frac{p}{2})^{(M-1)r}$, conditional only on $\xi$ & nothing further\\ \hline
\end{tabular}
\end{gridtable}
Two auxiliary quantities are used. The random variable $A_{i,n}$ records whether the routing entry on $e_i$ accepts the tested point $y_n(t)$, and $Z$ counts how many tested points are accepted by their routing entries. Once all the routing entries are fixed, a tested point with $A_{i,n}=0$ has already failed its local test, and each of the remaining $Z$ local tests fails precisely when the value of its terminal entry is zero. Acceptance alone, $A_{i,n}=1$, does not imply membership in $B$; success of the local test, $Q_i(y_n(t))=1$, does, by lemma~\ref{lem:route}, while failure of the local test does not imply exclusion from $B$.

\begin{figure}[htbp]
\centering
\begin{tikzpicture}[>=Stealth,font=\small,
  vertex/.style={circle,draw=black!60,fill=white,inner sep=2pt}]
\node[vertex] (Ui) at (0,3.2) {$U_i$};
\draw[black!40] (Ui)--(-2.2,.9)--(2.2,.9)--cycle;
\draw[leftbranch,line width=1.1pt] ($(Ui.south)+(-.06,0)$)--(-.41,2.0)--(.09,.9);
\draw[rightbranch,line width=1.1pt,dash pattern=on 3.5pt off 2pt] ($(Ui.south)+(.06,0)$)--(-.29,2.0)--(.21,.9);
\fill[black] (.15,.9) circle (.07);
\node[anchor=south west] at (.25,.9) {leaf $L$};
\node[anchor=west,align=left,font=\footnotesize] at (2.4,2.4) {$y$ (solid) and $y'$ (dashed) may follow\\the same route to the same leaf};
\node[anchor=east] at (-2.6,-.45) {$T_L$};
\foreach \k in {0,...,11}{\draw[black!60] ({-2.4+\k*.4},-.7) rectangle ({-2.0+\k*.4},-.2);}
\path[fill=leftbranch!30,draw=black!60] (-1.2,-.7) rectangle (-.8,-.2);
\path[fill=rightbranch!30,draw=black!60] (1.2,-.7) rectangle (1.6,-.2);
\draw[->,leftbranch,thick] (.15,.85)--(-1.0,-.18);
\draw[->,rightbranch,thick,dash pattern=on 3.5pt off 2pt] (.15,.85)--(1.4,-.18);
\node[below,font=\footnotesize,text=leftbranch] at (-1.0,-.75) {$J_{b(L)}(y)$};
\node[below,font=\footnotesize,text=rightbranch] at (1.4,-.75) {$J_{b(L)}(y')$};
\end{tikzpicture}
\caption{Two tested points $y$ and $y'$ attached to the same child $U_i$ may reach the same leaf $L$. Since paths in a tree are unique, they then follow the same edges from $U_i$ to $L$, although at each vertex on the way they read different routing entries, and they read the same terminal table $T_L$. They nevertheless read different entries of it: they have distinct cell labels in the grid of $e_i$ by lemma~\ref{lem:separation}, and hence in the finer grid of $L$.}
\label{fig:terminal-entries}
\end{figure}
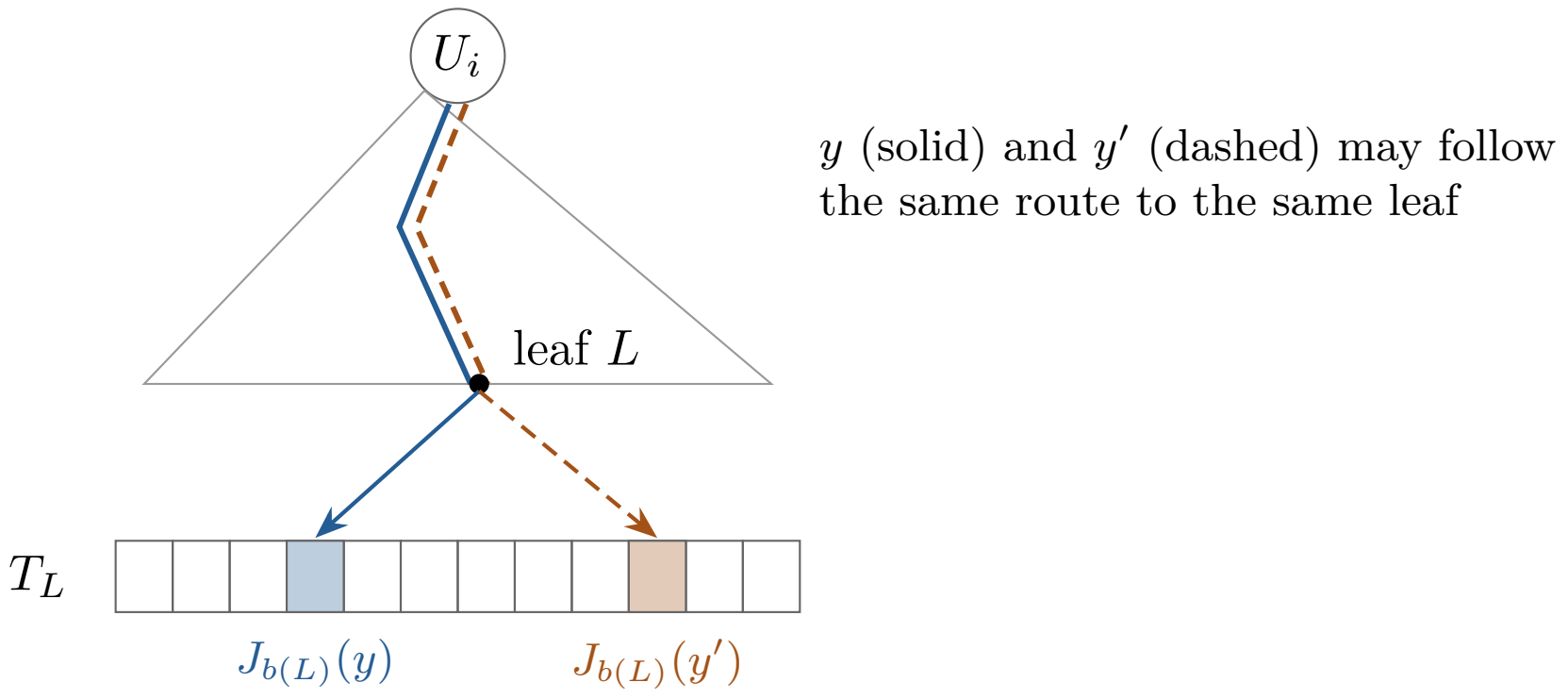

It is important to distinguish a table from an entry. A terminal entry is identified by a pair
\[
(\text{leaf table},\ \text{cell label}),
\]
and two tested points may read the same table while reading different entries of it; see figure~\ref{fig:terminal-entries}.

For two tested points, the calculation can be written out completely. The pair of their routing entries on the relevant edges takes the values $(0,0)$, $(0,1)$, $(1,0)$, $(1,1)$, each with probability $\frac{1}{4}$, and the corresponding conditional probabilities that both local tests fail are $1$, $1-p$, $1-p$ and $(1-p)^2$. Averaging gives
\[
\frac{1}{4}+\frac{1}{2}(1-p)+\frac{1}{4}(1-p)^2=\Bigl(1-\frac{p}{2}\Bigr)^2.
\]
The general calculation repeats this argument for $(M-1)r$ tested points.

\begin{proof}[Proof of lemma~\ref{lem:failure}]
For each $1\le i<M$ and each $n\in W_{e_i}$, consider the tested point $y_n(t)=x+tq^n$ and define
\[
A_{i,n}=S_{e_i}(J_{b_{e_i}}(y_n(t))).
\]
Thus $A_{i,n}$ records whether the routing entry on $e_i$ accepts this point: it equals one if the point is accepted and zero if it is rejected. For a fixed $i$, lemma~\ref{lem:separation} says that the cell labels of the tested points in the grid of $e_i$ are different from one another and from the cell label of $x$, so none of these entries was fixed by $\xi$. Since different values of $i$ use different routing tables, the $A_{i,n}$ are independent Bernoulli-$\frac{1}{2}$ random variables under the conditioning $\xi$.

Now reveal every remaining routing entry. For each tested point this determines a leaf $L_{i,n}=L_i(y_n(t))$ and a terminal entry, identified by the pair
\[
(L_{i,n},J_{b(L_{i,n})}(y_n(t))).
\]
These pairs are all distinct, for three reasons.
\begin{enumerate}
\item Tested points attached to different children of $U$ end in different subtrees, hence in different leaf tables.
\item Tested points attached to the same child but ending in different leaves also read different tables.
\item If two tested points attached to the same child reach the same leaf $L$, then $b(L)\ge b_{e_i}$. They have distinct cell labels in the grid of $e_i$, and lemma~\ref{lem:separation} preserves that distinction in the grid of $L$, so they read different entries of the same table.
\end{enumerate}
The third case is exactly where the geometry of section~\ref{sec:geometry} enters the probability calculation.

Consequently, after all the routing entries have been revealed, the terminal entries read by the tested points are independent Bernoulli-$p$ variables, none of them revealed by $\xi$. With
\[
Z=\sum_{i=1}^{M-1}\sum_{n\in W_{e_i}}A_{i,n},
\]
the conditional probability that every local test fails is $(1-p)^Z$.

Finally, average this expression over the remaining routing entries, still holding $\xi$ fixed. For one Bernoulli-$\frac{1}{2}$ random variable $A$, the expectation of $(1-p)^A$ is
\[
\frac{1}{2}+\frac{1}{2}(1-p)=1-\frac{p}{2}.
\]
Since the $(M-1)r$ random variables $A_{i,n}$ are independent Bernoulli-$\frac{1}{2}$ variables under the conditioning $\xi$, the expectation of the product $(1-p)^Z=\prod_{i,n}(1-p)^{A_{i,n}}$ is the product of these expectations:
\[
{\mathbb P}(\text{all the local tests fail}\mid\xi)
={\mathbb E}\bigl[(1-p)^Z\mid\xi\bigr]
=\prod_{i=1}^{M-1}\prod_{n\in W_{e_i}}{\mathbb E}\bigl[(1-p)^{A_{i,n}}\mid\xi\bigr]
=\Bigl(1-\frac{p}{2}\Bigr)^{(M-1)r}.
\]
The first equality averages the conditional probability $(1-p)^Z$, obtained after fixing all the routing entries, over the routing entries not fixed by $\xi$; the second uses the conditional independence of the $A_{i,n}$. This proves~\eqref{eq:exact-failure}.
\end{proof}

The calculation separates two uses of randomness: the fair routing entries make many independent opportunities available, while the terminal Bernoulli variables keep the expected density small and decide whether those opportunities succeed. Neither the number of branches nor the number of tested points changes the identity ${\mathbb E}\rho(B)=p$.

\section{How many dilation values must be tested?}\label{sec:scales}

Lemma~\ref{lem:failure} controls one fixed dilation $t$. A bound valid for each fixed $t$ does not by itself control the event that a failure occurs at some $t\in[1,2]$, since there are uncountably many dilations. The mechanism which overcomes this is the following.
\begin{quote}
As $t$ varies, the relevant cell labels change only at finitely many boundary-crossing times. Between consecutive crossing times, every local test has the same outcome, for every realization of the tables.
\end{quote}
This is an exact finite reduction: we do not approximate $t$ by a nearby value and estimate an approximation error. Counting the finitely many ways in which a dilation can change cell membership, and then applying the union bound, is the approach of Chleb\'ik~\cite[section~5]{Chlebik2015} and Kolountzakis--Papageorgiou~\cite[section~3.1]{KolountzakisPapageorgiou2025}; Solymosi's report~\cite[section~3]{Solymosi2011} likewise studies finite choices of the dilation. In the OpenAI proof, the count only involves the grids consulted by the local tests.

Throughout the section we fix, as in lemma~\ref{lem:failure}, a stable base point $x$ and a specification $\xi$ of the routing entries at $x$ for which the route of $x$ has a first default vertex $U$. Then $U$, its height $h$, the length $r=r_h$, the edges $e_i=(U,i)$ with $1\le i<M$, and their index intervals $W_{e_i}$ are fixed, while the entries not revealed by $\xi$ are still random. Recall that the tested points are the points $y_n(t)=x+tq^n$ with $n\in W_{e_i}$ and $1\le i<M$, see~\eqref{eq:tested-points}, and that the local tests $Q_i$ are defined in~\eqref{eq:local-test}.

\subsection{Which grid must be watched}\label{subsec:local-grid}

The local test $Q_i(z)$ reads only three kinds of tables: the routing table on $e_i$, the routing tables on the edges below $U_i$, and the terminal tables at the leaves below $U_i$. Recall from~\eqref{eq:bstar} that $b_{e_i}^*$ is the largest index assigned to $e_i$ or to an edge below its child $U_i$; here ``below $U_i$'' refers to all the edges of the subtree below $U_i$, at every depth, and not only to the edges leaving $U_i$. Each of the tables just listed uses a grid indexed by at most $b_{e_i}^*$: the routing table on an edge $f$ uses the grid indexed by $b_f$, and the terminal table at a leaf $L$ uses the grid indexed by $b(L)=b_f$, where $f$ is the edge entering $L$, which is $e_i$ or an edge below $U_i$. Since the grids are nested, the grid $\Gamma_{b_{e_i}^*}$ refines all of them, and by~\eqref{eq:nesting} the cell label $J_{b_{e_i}^*}(z)$ determines every coarser label that $Q_i$ might read. We call $\Gamma_{b_{e_i}^*}$ the \emph{local grid} of $e_i$, and of the tested points attached to $e_i$. Once the tables are fixed, the value $Q_i(z)$ is therefore determined by the single cell label $J_{b_{e_i}^*}(z)$. This is the reason the local tests were introduced: they give sufficient evidence of membership in $B$ while consulting only a specified part of the tree, and this restricted dependence is what makes a local count possible.

\subsection{Counting the crossings}\label{subsec:crossings}

Fix an index $n\in W_{e_i}$ of a tested point. We call $t\in[1,2]$ a \emph{crossing time} of this tested point if $x+tq^n\in\Gamma_{b_{e_i}^*}$, that is, if at the dilation $t$ the point lies on a boundary of its local grid. As $t$ goes from $1$ to $2$, the point $x+tq^n$ moves strictly to the right through the closed interval $[x+q^n,x+2q^n]$, of length $q^n$, so distinct boundaries in this interval correspond to distinct crossing times. A closed interval of length $\ell$ contains at most $1+\ell N$ points of $\frac{1}{N}{\mathbb Z}$, and the cells of the local grid have width $1/N_{b_{e_i}^*}$, so the number of crossing times of this tested point is at most
\[
1+\frac{\text{distance travelled}}{\text{cell width}}=1+q^nN_{b_{e_i}^*}.
\]
Only now do we substitute the grid estimate~\eqref{eq:grid-size}:
\begin{equation}\label{eq:crossings}
1+q^nN_{b_{e_i}^*}
\le1+2cq^{n-b_{e_i}^*}
\le1+2cq^{-2r}.
\end{equation}
For the last inequality we use $n\ge a_{e_i}$ and the span estimate $b_{e_i}^*-a_{e_i}+1\le2r$ of lemma~\ref{lem:span}, which applies because $e_i$ leaves the vertex $U$ of height $h$ and $r=r_h$; it controls the difference $b_{e_i}^*-n$. This is where the construction of section~\ref{subsec:indices} pays off: the grid indexed by $b_{e_i}^*$ may be extremely fine, but the displacement $q^n$ of the point is correspondingly small. We stress that $q^{-2r}$ bounds a number of boundary crossings; it is not the resolution of the grid, which is about $q^{-b_{e_i}^*}$.

All the tested points use the same dilation parameter $t$. We therefore merge their crossing times into one list on the single interval $[1,2]$, and their numbers add: with $(M-1)r$ tested points there are at most $(M-1)r(1+2cq^{-2r})$ crossing times in total. We do not multiply the numbers of cell configurations of the separate tested points.

\subsection{Representatives}\label{subsec:representatives}

We now replace the interval $[1,2]$ of dilations by a finite set. Let $k$ be the number of distinct crossing times of the tested points which lie in the open interval $(1,2)$, and let
\begin{equation}\label{eq:listed-times}
1=t_0<t_1<\cdots<t_{k+1}=2
\end{equation}
be these crossing times together with the endpoints $1$ and $2$; we call them the \emph{listed times}. The set of \emph{representatives} consists of the listed times and of one midpoint between each pair of consecutive listed times:
\begin{equation}\label{eq:representatives}
{\mathcal T}_{x,\xi}=\{t_0,t_1,\ldots,t_{k+1}\}\cup\Bigl\{\frac{t_{l-1}+t_l}{2}:1\le l\le k+1\Bigr\}.
\end{equation}
The representatives are useful because of the following property:
\begin{equation}\label{eq:representative-property}
\begin{gathered}
\text{for every }t\in[1,2]\text{ there is }\tau\in{\mathcal T}_{x,\xi}\text{ such that}\\
J_{b_{e_i}^*}(x+tq^n)=J_{b_{e_i}^*}(x+\tau q^n)\qquad\text{for all }i<M\text{ and }n\in W_{e_i}.
\end{gathered}
\end{equation}
Indeed, if $t$ is a listed time, we take $\tau=t$. Otherwise $t$ lies in an open interval $(t_{l-1},t_l)$ between consecutive listed times, and we take for $\tau$ the midpoint of this interval. No crossing time lies in $(t_{l-1},t_l)$, so as $s$ runs through this interval, no tested point $x+sq^n$ lies on a boundary of its local grid; since the point moves continuously, it stays in one cell, and its label is the same at $s=t$ as at $s=\tau$. By section~\ref{subsec:local-grid}, once the tables are fixed each local test is determined by these labels, so in every realization of the tables all the local tests have the same outcomes at $t$ as at $\tau$.

The listed times alone would also suffice: since the points move to the right as $t$ increases, and a point lying on a boundary belongs to the cell on its right, every label is constant on each half-open interval $[t_{l-1},t_l)$, so $t_{l-1}$ could represent that whole interval. We keep the midpoints because they make the argument above independent of this convention, at the cost of at most doubling the number of representatives.

We now count the representatives. There are $k+2$ listed times and $k+1$ midpoints, so $\#{\mathcal T}_{x,\xi}=2k+3$. The $k$ listed times in $(1,2)$ are crossing times of tested points, so $k$ is at most the total number of crossing times of all the tested points, which is at most $(M-1)r(1+2cq^{-2r})$ by section~\ref{subsec:crossings}. Therefore
\begin{equation}\label{eq:D-def}
\#{\mathcal T}_{x,\xi}\le D_{M,q}(r)=3+2\underbrace{(M-1)r}_{\text{number of tested points}}\,\underbrace{(1+2cq^{-2r})}_{\text{crossing bound per tested point}}.
\end{equation}

Finally, ${\mathcal T}_{x,\xi}$ is not random. Its elements are computed from the base point $x$, the deterministic grids and index intervals, and the vertex $U$; this vertex is determined by $\xi$, and no other table entry is involved. This is what allows us to use lemma~\ref{lem:failure} at the representatives: that lemma concerns one dilation fixed in advance, before the unrevealed entries are averaged, and once $x$ and $\xi$ are fixed, each $\tau\in{\mathcal T}_{x,\xi}$ is such a dilation. If the representatives depended on unrevealed entries, lemma~\ref{lem:failure} could not be applied at them.

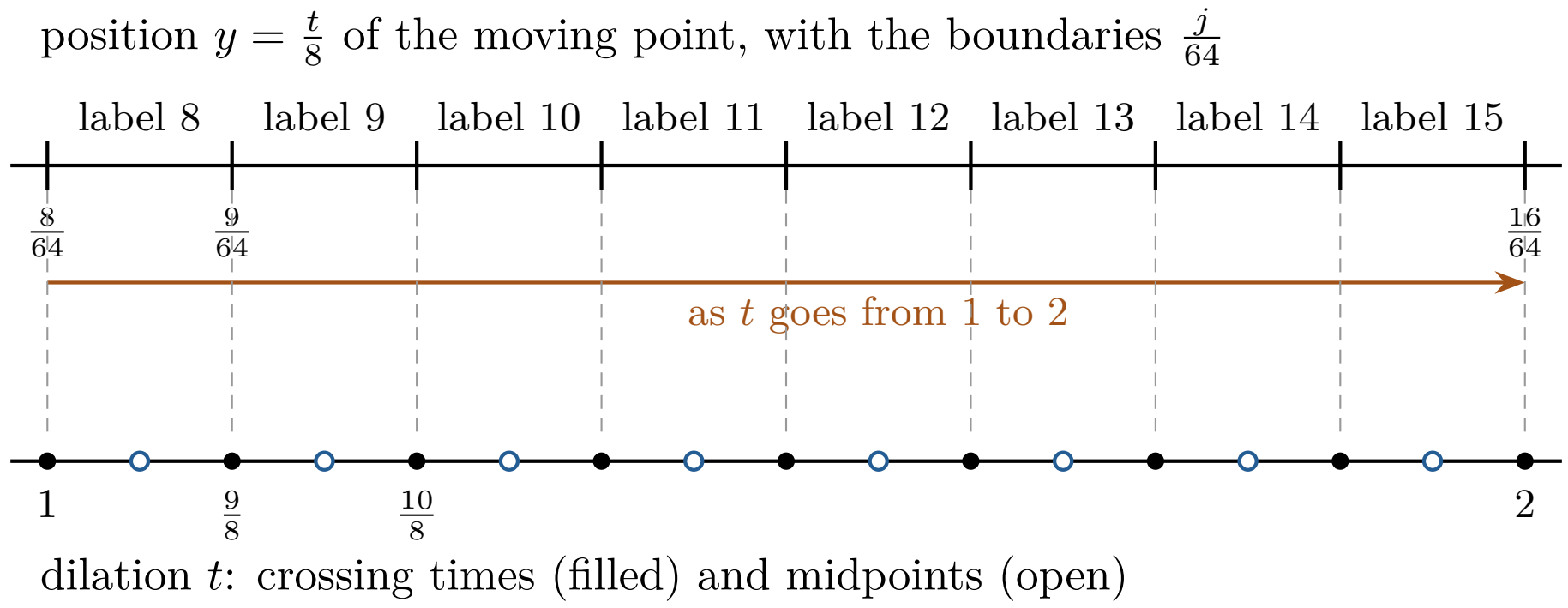
\begin{figure}[htbp]
\centering
\begin{tikzpicture}[x=1cm,y=1cm,>=Stealth,font=\small,scale=.95]
\node[anchor=west] at (-.2,2.85) {position $y=\frac{t}{8}$ of the moving point, with the boundaries $\frac{j}{64}$};
\draw[thick] (-.3,1.8)--(12.3,1.8);
\foreach \j in {0,...,8}{\draw[thick] ({1.5*\j},1.6)--({1.5*\j},2.0);}
\foreach \j/\lab in {0/8,1/9,2/10,3/11,4/12,5/13,6/14,7/15}{\node[font=\footnotesize] at ({1.5*\j+.75},2.2) {label $\lab$};}
\node[below,font=\footnotesize] at (0,1.6) {$\frac{8}{64}$};
\node[below,font=\footnotesize] at (1.5,1.6) {$\frac{9}{64}$};
\node[below,font=\footnotesize] at (12,1.6) {$\frac{16}{64}$};
\draw[->,rightbranch,thick] (0,.85)--(12,.85);
\node[below,font=\footnotesize,text=rightbranch,fill=white,inner sep=1.5pt] at (6.75,.8) {as $t$ goes from $1$ to $2$};
\foreach \j in {0,...,8}{\draw[densely dashed,black!40] ({1.5*\j},1.6)--({1.5*\j},-.4);}
\draw[thick] (-.3,-.6)--(12.3,-.6);
\foreach \j in {0,...,8}{\fill ({1.5*\j},-.6) circle (.07);}
\foreach \j in {0,...,7}{\draw[leftbranch,fill=white,thick] ({1.5*\j+.75},-.6) circle (.07);}
\node[below,font=\footnotesize] at (0,-.7) {$1$};
\node[below,font=\footnotesize] at (1.5,-.7) {$\frac{9}{8}$};
\node[below,font=\footnotesize] at (3,-.7) {$\frac{10}{8}$};
\node[below,font=\footnotesize] at (12,-.7) {$2$};
\node[anchor=west] at (-.2,-1.55) {dilation $t$: crossing times (filled) and midpoints (open)};
\end{tikzpicture}
\caption{The illustrative case $q=\frac{1}{2}$, $x=0$, $n=3$, $N_3=64$. Above, the point $\frac{t}{8}$ crosses the boundaries $\frac{j}{64}$, and its cell label increases from $8$ to $15$ before reaching the boundary $\frac{16}{64}$ at $t=2$. Below, the corresponding crossing times $t=\frac{j}{8}$ on the dilation axis, joined to the boundaries by vertical guides. The local test is constant on each open interval between crossing times, and the filled and open dots are the $17$ representatives. With several tested points, all their crossing times are merged on the same axis.}
\label{fig:representatives}
\end{figure}

For example, take $q=\frac{1}{2}$, so that $c=8$, and consider, purely as an illustration, a single tested point with $x=0$ and $n=b^*=3$, for which $N_3=64$. The moving point is $tq^3=\frac{t}{8}$, and its cell label is
\[
J_3\Bigl(\frac{t}{8}\Bigr)=\lfloor8t\rfloor.
\]
The crossing times are $t=1,\frac{9}{8},\frac{10}{8},\ldots,2$, nine in all, which is exactly the bound $1+q^3N_3=1+8$. At $t=\frac{9}{8}$, for example, the point reaches the boundary $\frac{9}{64}$ and its label changes from $8$ to $9$. Here $k=7$, and the $9$ crossing times together with $8$ midpoints give $17=2k+3$ representatives; see figure~\ref{fig:representatives}.

\subsection{The estimate for all dilations}

Recall from~\eqref{eq:N-def} that ${\mathcal N}=\bigcup_eW_e$ is the finite set of all the indices assigned to the edges of the tree; it contains the index intervals $W_{e_i}$ of the tested points. The event that we have to control is the \emph{failure event} at $x$,
\begin{equation}\label{eq:failure-event}
F_x=\bigl\{\text{some }t\in[1,2]\text{ satisfies }x+tq^n\notin B\text{ for every }n\in{\mathcal N}\bigr\},
\end{equation}
that is, the event that for some dilation all the points $x+tq^n$ with $n\in{\mathcal N}$ miss $B$. At such a dilation the tested points miss $B$ as well, so by lemma~\ref{lem:route} every local test fails there; the converse is not needed. Combined with the finite reduction of section~\ref{subsec:representatives}, this gives the following estimate.

\begin{lemma}\label{lem:all-scales}
Let $x$ be a stable base point, and let $\xi$ be a specification of the routing entries at $x$ for which the route of $x$ has a first default vertex $U$, of height $h$; put $r=r_h$. Then
\begin{equation}\label{eq:all-scales}
{\mathbb P}(F_x\mid\xi)\le D_{M,q}(r)e^{-\frac{p(M-1)r}{2}},
\end{equation}
where $D_{M,q}(r)=3+2(M-1)r(1+2cq^{-2r})$ is the bound~\eqref{eq:D-def} on the number of representatives.
\end{lemma}

\begin{proof}
All the probabilities in this proof are conditional on $\xi$, so only the realizations of the tables in which the routing entries at the cell labels of $x$ take the values given by $\xi$ are considered. For $\tau\in[1,2]$, let
\[
{\mathcal E}_\tau=\bigl\{Q_i(y_n(\tau))=0\text{ for every }i<M\text{ and every }n\in W_{e_i}\bigr\}
\]
be the event that all the local tests fail at the dilation $\tau$.

\emph{Step 1: reduction to the representatives.} Fix a realization of the tables for which some $t\in[1,2]$ makes every point $x+tq^n$ with $n\in{\mathcal N}$ miss $B$, that is, a realization in $F_x$. Since $W_{e_i}\subset{\mathcal N}$, every tested point $y_n(t)$ misses $B$, so by lemma~\ref{lem:route} all the local tests fail at $t$; the realization lies in ${\mathcal E}_t$. By~\eqref{eq:representative-property}, there is $\tau\in{\mathcal T}_{x,\xi}$ at which every local test has the same outcome as at $t$, so the realization also lies in ${\mathcal E}_\tau$. We have shown that
\begin{equation}\label{eq:reduction}
F_x\subset\bigcup_{\tau\in{\mathcal T}_{x,\xi}}{\mathcal E}_\tau.
\end{equation}

\emph{Step 2: the union bound.} The set ${\mathcal T}_{x,\xi}$ is finite, and, as explained at the end of section~\ref{subsec:representatives}, it is determined by $x$ and $\xi$, so it is the same for all the realizations under consideration. The union bound, that is, the subadditivity of the conditional probability ${\mathbb P}(\,\cdot\mid\xi)$, applied to the finite union in~\eqref{eq:reduction}, gives
\[
{\mathbb P}(F_x\mid\xi)\le\sum_{\tau\in{\mathcal T}_{x,\xi}}{\mathbb P}({\mathcal E}_\tau\mid\xi).
\]
No independence between the events ${\mathcal E}_\tau$ for different $\tau$ is needed here.

\emph{Step 3: one representative at a time.} For each fixed $\tau\in{\mathcal T}_{x,\xi}$, lemma~\ref{lem:failure}, applied with $t=\tau$, gives
\[
{\mathbb P}({\mathcal E}_\tau\mid\xi)=\Bigl(1-\frac{p}{2}\Bigr)^{(M-1)r}.
\]
All the terms of the sum are therefore equal, and by the count~\eqref{eq:D-def},
\[
{\mathbb P}(F_x\mid\xi)\le\#{\mathcal T}_{x,\xi}\,\Bigl(1-\frac{p}{2}\Bigr)^{(M-1)r}\le D_{M,q}(r)\Bigl(1-\frac{p}{2}\Bigr)^{(M-1)r}.
\]

\emph{Step 4: the exponential form.} The elementary inequality $1-u\le e^{-u}$, valid for every real $u$, gives $0<1-\frac{p}{2}\le e^{-\frac{p}{2}}$, since $p<1$. Raising both sides to the power $(M-1)r$ gives $\bigl(1-\frac{p}{2}\bigr)^{(M-1)r}\le e^{-\frac{p(M-1)r}{2}}$, and~\eqref{eq:all-scales} follows.
\end{proof}

The distinction between the local grid of $e_i$ and the finest grid of the whole tree is decisive here. Counting crossings in the finest grid would produce a factor which grows with the size of the whole tree, while in the local grid the exponent in~\eqref{eq:crossings} is at most $2r$, however far down the index list the other branches extend. To make this quantitative, let $r\ge1$ be an integer. Then $1\le r\le rq^{-2r}$, so
\[
D_{M,q}(r)=3+2(M-1)r+4c(M-1)rq^{-2r}\le C_{M,q}\,r\,q^{-2r},\qquad\text{where}\quad C_{M,q}=3+(M-1)(2+4c)
\]
and $c=\frac{4}{1-q}$ is the constant of section~\ref{subsec:grids}. Writing $q^{-2r}=e^{2r\log\frac{1}{q}}$, we obtain
\begin{equation}\label{eq:rate}
D_{M,q}(r)e^{-\frac{p(M-1)r}{2}}\le C_{M,q}\,r\,
\exp\Bigl(-\Bigl[\frac{p(M-1)}{2}-2\log\frac{1}{q}\Bigr]r\Bigr)\qquad\text{for every integer }r\ge1.
\end{equation}
This explains the choice of parameters in section~\ref{sec:finish}: first choose $M$ so large that the bracket is positive, and then let $r$ grow.

The gain over the single grid of section~\ref{sec:obstacle} lies in this bracket. In both constructions, a union bound multiplies a failure probability which decays exponentially in the number of tests by a number of cell configurations, or of representative dilations, which grows exponentially, and the question is which exponent is larger. With one grid, $N$ tested points in distinct cells all fail with probability at most $e^{-pN}$. By~\eqref{eq:gap-obstruction}, however, the smallest gap between $N$ powers of $q$ is at most $\frac{q^{N-2}}{1-q}$ times their spread, so cells which separate them are correspondingly small, and the number of cell configurations met as the dilation varies grows at least like $q^{-N}$, up to factors independent of $N$. The union bound is then of the order of $\exp\bigl((\log\frac{1}{q}-p)N\bigr)$, which does not tend to zero when $p<\log\frac{1}{q}$. Since $p$ is the density of the set we remove, it must be allowed to be small, and the single grid has no parameter which changes the ratio between the two exponents. In the tree, the failure exponent is $\frac{p(M-1)}{2}r$, while the logarithm of the number of representative dilations is at most $2r\log\frac{1}{q}+\log r+O_{M,q}(1)$. Increasing the branching number $M$ multiplies the number of tests, and hence the failure exponent, but changes neither the coefficient $2\log\frac{1}{q}$ of $r$ in the logarithm of the count nor, by lemma~\ref{lem:density}, the expected density $p$ of $B$. Choosing $M$ with $\frac{p(M-1)}{2}>2\log\frac{1}{q}$ therefore makes the decay win, for every $p\in(0,1)$. This is the gain which a single grid did not provide.

The same comparison, between the number of separated points available and the logarithm of the number of representative dilations, drives the criterion for finite subsets of Kolountzakis and Chleb\'ik and its use in the two-fold sumset paper~\cite{IKMRY2026Sumset}; section~\ref{sec:comparison} explains this parallel in detail.

\section{Descent of the unresolved probability}\label{sec:descent}

Sections~\ref{sec:independence} and~\ref{sec:scales} fixed the specification $\xi$ of the routing entries at the base point. We now keep the stable base point $x$ fixed but allow its routing entries to vary again, and we combine the estimates of section~\ref{sec:scales} according to where the route of $x$ first makes a default choice. A convenient way to organize this combination comes from the viewpoint of Burstein, Iosevich, and Krause~\cite{BursteinIosevichKrause2026}: one retains a single nonnegative quantity throughout an exploration, and uses its decreases to weight the contributions made where the exploration stops. Here this bookkeeping combines the local estimates without losing a factor proportional to the depth of the tree.

\subsection{The probability that no default has yet occurred}

Recall from section~\ref{sec:routing} that the route of $x$ makes one routing decision at each vertex it visits, namely the choice of a child, and that this decision is a default choice when all the $M-1$ routing entries read by $x$ at that vertex are zero. The route visits the vertices of depths $0,1,\ldots,d-1$ before it reaches a leaf, so it makes exactly $d$ routing decisions, the $(j+1)$-st being made at depth $j$. Fix a stable base point $x$, and set
\begin{equation}\label{eq:theta}
\theta=2^{-(M-1)},
\end{equation}
the probability that $M-1$ independent fair routing entries are all zero. For $0\le j\le d$, let
\begin{equation}\label{eq:Omega}
\Omega_j=\{\text{the first $j$ routing decisions of $x$ are all nondefault}\},
\end{equation}
so that $\Omega_0$ is the whole probability space and $\Omega_0\supset\Omega_1\supset\cdots\supset\Omega_d$, and define
\begin{equation}\label{eq:survival-potential}
u_j={\mathbb P}(\Omega_j).
\end{equation}
On $\Omega_j$, the route has not yet reached a vertex at which the local estimate of section~\ref{sec:scales} can be applied (this is the reason to be called ``unresolved'').

Let $0\le j<d$, and condition on the routing entries read by $x$ during its first $j$ routing decisions, on the event $\Omega_j$. These entries determine the vertex $P$ of depth $j$ reached by the route of $x$. The decision at $P$ is made with the $M-1$ entries $S_{(P,1)}(J_{b_{(P,1)}}(x)),\ldots,S_{(P,M-1)}(J_{b_{(P,M-1)}}(x))$, which belong to tables not read at the earlier vertices; hence they are independent of the conditioning, and they are all zero with probability $\theta$. Averaging over the conditioning gives
\[
u_{j+1}={\mathbb P}(\Omega_{j+1})=(1-\theta)\,{\mathbb P}(\Omega_j)=(1-\theta)u_j.
\]
Since $u_0={\mathbb P}(\Omega_0)=1$, induction on $j$ gives
\begin{equation}\label{eq:u-formula}
u_j=(1-\theta)^j\qquad(0\le j\le d);
\end{equation}
for $j=d$ this is~\eqref{eq:no-default}. We call $u_j$ a \emph{potential}: it is nonnegative and decreases with $j$, and its decrease
\begin{equation}\label{eq:potential-decrease}
u_j-u_{j+1}={\mathbb P}(\Omega_j\setminus\Omega_{j+1})=\theta(1-\theta)^j
\end{equation}
is exactly the probability that the first default occurs at depth $j$, the root having depth zero. A vertex at depth $j$ has height $d-j$.

\subsection{Combining the local estimates}

For $1\le h\le d$, let
\begin{equation}\label{eq:local-error}
\varepsilon_h=D_{M,q}(r_h)e^{-\frac{p(M-1)r_h}{2}},
\end{equation}
which is the bound of lemma~\ref{lem:all-scales} when the first default vertex has height $h$, so that $r=r_h$; we call $\varepsilon_h$ the \emph{local error} at height $h$. Recall the failure event $F_x$ of~\eqref{eq:failure-event}: some dilation $t\in[1,2]$ makes all the points $x+tq^n$ with $n\in{\mathcal N}$ miss $B$. By lemma~\ref{lem:all-scales}, every specification $\xi$ whose first default vertex has height $h$ satisfies ${\mathbb P}(F_x\mid\xi)\le\varepsilon_h$. The possible routes of $x$ fall into the following mutually exclusive cases.
\begin{gridtable}
\begin{tabular}{|>{\raggedright\arraybackslash}p{4.6cm}|>{\raggedright\arraybackslash}p{3.0cm}|>{\raggedright\arraybackslash}p{4.4cm}|}
\hline
\emph{what happens to the route} & \emph{probability} & \emph{available bound on ${\mathbb P}(F_x\mid\cdot)$}\\ \hline
first default at depth $j$, $0\le j<d$ & $u_j-u_{j+1}$ & $\varepsilon_{d-j}$\\ \hline
no default in all $d$ decisions & $u_d$ & $1$\\ \hline
\end{tabular}
\end{gridtable}
The subscript $d-j$ is the height of the first default vertex.

The law of total probability, applied to the cases in the table, gives
\begin{equation}\label{eq:weighted-descent}
{\mathbb P}(F_x)
\le u_d+\sum_{j=0}^{d-1}(u_j-u_{j+1})\varepsilon_{d-j}.
\end{equation}
This requires no independence between the stopping depth and the failure event. To see this, recall from~\eqref{eq:Omega} that $\Omega_j$ is the event that the first $j$ routing decisions of $x$ are all nondefault. Thus, for $0\le j<d$, the event $\Omega_j\setminus\Omega_{j+1}$ is the event that the first default occurs at depth $j$, while $\Omega_d$ is the event that no default occurs. The specification $\xi$ determines the route of $x$, so $\Omega_j\setminus\Omega_{j+1}$ is the union of the events $\{\xi\}$ over the set $\Xi_j$ of specifications whose first default occurs at depth $j$; here $\{\xi\}$ denotes the event~\eqref{eq:specification}, and ${\mathbb P}(\xi)$ its probability. Lemma~\ref{lem:all-scales} bounds ${\mathbb P}(F_x\mid\xi)$ by the same number $\varepsilon_{d-j}$ for every $\xi\in\Xi_j$, so
\[
{\mathbb P}\bigl(F_x\cap(\Omega_j\setminus\Omega_{j+1})\bigr)
=\sum_{\xi\in\Xi_j}{\mathbb P}(F_x\mid\xi)\,{\mathbb P}(\xi)
\le\varepsilon_{d-j}\sum_{\xi\in\Xi_j}{\mathbb P}(\xi)
=\varepsilon_{d-j}(u_j-u_{j+1}).
\]
The remaining specifications are those without a default, which make up $\Omega_d$, so they contribute at most ${\mathbb P}(\Omega_d)=u_d$.

If a default occurs, its first occurrence is unique, so the cases in the table are disjoint, and the decreases telescope:
\[
\sum_{j=0}^{d-1}(u_j-u_{j+1})=1-u_d.
\]
Hence if every local error satisfies $\varepsilon_h\le\varepsilon$ for some $0\le\varepsilon\le1$, then their total contribution is at most $\varepsilon(1-u_d)$, and~\eqref{eq:weighted-descent}, together with $u_d=(1-\theta)^d$ from~\eqref{eq:u-formula}, becomes
\begin{equation}\label{eq:descent-bound}
{\mathbb P}(F_x)\le u_d+\varepsilon(1-u_d)
=\varepsilon+(1-\varepsilon)(1-\theta)^d.
\end{equation}
The local errors are averaged with total weight at most one; they are not paid once for every possible stopping depth.

The estimate~\eqref{eq:descent-bound} is all that section~\ref{sec:finish} needs from this section. The next two subsections rewrite it as a recurrence and explain its connection with the work of Burstein, Iosevich, and Krause.

\FloatBarrier
\subsection{The same calculation as a recurrence}

Once the weighted bound is understood, it can be rewritten as a recurrence. For $\ell\ge0$, let $v_\ell$ be the bound obtained when $\ell$ routing decisions remain available, defined by
\[
v_0=1,\qquad v_{\ell+1}=\theta\varepsilon+(1-\theta)v_\ell.
\]
The interpretation is immediate: the next decision is a default with probability $\theta$, which contributes the allowance $\varepsilon$, or it is nondefault, leaving $\ell$ decisions available. The value $v_\ell$ is a numerical upper bound, not an actual conditional failure probability. Solving the recurrence gives
\[
v_\ell=\varepsilon+(1-\varepsilon)(1-\theta)^\ell,
\]
so that $v_d$ is the right side of~\eqref{eq:descent-bound}, and the excess above the local error floor satisfies
\[
v_{\ell+1}-\varepsilon=(1-\theta)(v_\ell-\varepsilon).
\]
Thus $v_\ell-\varepsilon$ is a geometrically decreasing nonnegative potential. The recurrence repackages the first-default calculation in the OpenAI proof and requires no new probabilistic hypothesis.

\subsection{The precise connection with Burstein--Iosevich--Krause}

The greedy proof of the selection theorem in~\cite[sections~2.2 and~2.4]{BursteinIosevichKrause2026} builds an approximating function step by step: whenever the current approximation violates one of the required inequalities, it is corrected. A nonnegative quantity, which there is a convex functional of the current approximation and which we again call a potential, decreases at each correction. Since it starts from a bounded value and remains nonnegative, its decrease bounds the total number of corrections, counted with suitable weights. Finally, each possible sequence of corrections is recorded by one of finitely many descriptions, and the probability estimate is obtained by summing over these descriptions. The same three ideas organize the argument of this section:
\begin{gridtable}
\begin{tabular}{|>{\raggedright\arraybackslash}p{6.0cm}|>{\raggedright\arraybackslash}p{6.4cm}|}
\hline
\emph{organizing idea in the cited work} & \emph{its realization here}\\ \hline
retain a nonnegative potential & the probability $u_j$ that the route of $x$ has made no default in its first $j$ decisions\\ \hline
account for contributions through its decreases & weight the first-default cases by $u_j-u_{j+1}$\\ \hline
sum a probability estimate over finite descriptions & encode the dilations by the representatives of section~\ref{subsec:representatives}, that is, by cell configurations in the local grids\\ \hline
\end{tabular}
\end{gridtable}
The objects involved are different. Here the potential $u_j$ is a probability rather than a convex functional. Nothing is corrected: $u_j$ decreases because each routing decision is a default with probability $\theta$, see~\eqref{eq:u-formula}, and its decreases~\eqref{eq:potential-decrease} are exactly the weights of the cases in~\eqref{eq:weighted-descent}. Finally, a representative dilation records the exact cell labels of the tested points at that dilation, whereas a description in the cited work records a sequence of corrections, whose sum gives the final approximation. The resemblance between the two arguments is therefore one of organization, not of the objects to which this organization is applied.

The two ingredients work together. The geometric span estimate $b_e^*-a_e+1\le2r_h$ of lemma~\ref{lem:span} makes each local error $\varepsilon_h$ small, and the potential argument combines those errors without a factor proportional to the depth $d$.

\section{Choose the parameters and remove the exceptional base points}\label{sec:finish}

We can now prove proposition~\ref{prop:periodic}. Recall the goal: for a given $p\in(0,1)$, we want one open $1$-periodic set $H$ with density $\rho(H)\le6p$ which intersects every normalized copy $x+tG_q$, $x\in{\mathbb R}$, $t\in[1,2]$.

\subsection{The parameters}

Each parameter has a single task:
\begin{gridtable}
\begin{tabular}{|>{\raggedright\arraybackslash}p{2.6cm}|>{\raggedright\arraybackslash}p{9.8cm}|}
\hline
\emph{parameter} & \emph{what its choice accomplishes}\\ \hline
$M$ & the exponential decrease of the failure probability outweighs the growth of the number of cell configurations met as the dilation varies\\ \hline
$d$ & the probability of seeing no default along the whole route is small\\ \hline
$g$ & the set of unstable base points is small\\ \hline
$r_0$ & the local failure bound $\varepsilon_h$ is small at every height\\ \hline
\end{tabular}
\end{gridtable}
The order matters. Increasing $M$ makes defaults rarer, since $\theta=2^{-(M-1)}$, so $d$ must be chosen after $M$; once $M$ and $d$ are fixed, the number $\kappa$ of edges is fixed, and $g$ can be chosen in terms of $\kappa$.

Fix $p\in(0,1)$. Choose $M$ so large that
\begin{equation}\label{eq:alpha}
\alpha=\frac{p(M-1)}{2}>2\log\frac{1}{q}.
\end{equation}
With this $M$ fixed, choose $d$ so large that
\begin{equation}\label{eq:depth-choice}
(1-2^{-(M-1)})^d<p.
\end{equation}
The number $\kappa$ of edges is now fixed. Choose $g$ so large that
\begin{equation}\label{eq:gap-choice}
4\kappa cq^{g+1}<p.
\end{equation}
Finally, choose $r_0$ so large that
\begin{equation}\label{eq:window-choice}
D_{M,q}(r)e^{-\alpha r}<p
\qquad\text{for every integer }r\ge r_0.
\end{equation}
Here $D_{M,q}(r)$ is the bound~\eqref{eq:D-def} on the number of representative dilations. The different heights use different lengths $r_h$, and the requirement for every $r\ge r_0$ lets one choice of $r_0$ cover all of them. Such a choice exists. Indeed, since $\alpha=\frac{p(M-1)}{2}$ by~\eqref{eq:alpha}, the estimate~\eqref{eq:rate} reads
\[
D_{M,q}(r)e^{-\alpha r}\le C_{M,q}\,r\,
e^{-(\alpha-2\log\frac{1}{q})r}\qquad\text{for every integer }r\ge1.
\]
Here $M$ and $q$ are fixed, so $C_{M,q}$ is a constant, and the exponent $\alpha-2\log\frac{1}{q}$ is positive by~\eqref{eq:alpha}. Hence the right side tends to $0$ as the integer $r$ tends to infinity, so there is an integer $r_0$ such that it is less than $p$ for every integer $r\ge r_0$, and~\eqref{eq:window-choice} holds for this $r_0$.
Construct the intervals and grids with these parameters; all the $r_h$ are at least $r_0$.

Let $x$ be a stable base point, that is, $x\in{\mathcal S}$ in the notation of section~\ref{subsec:stable}. By~\eqref{eq:u-formula} and~\eqref{eq:depth-choice}, $u_d=(1-\theta)^d<p$. By~\eqref{eq:local-error} and~\eqref{eq:alpha}, the local error at height $h$ is $\varepsilon_h=D_{M,q}(r_h)e^{-\alpha r_h}$, which is less than $p$ by~\eqref{eq:window-choice}, since $r_h\ge r_0$. Applying~\eqref{eq:descent-bound} with $\varepsilon=p$ gives
\begin{equation}\label{eq:pointwise-miss}
{\mathbb P}(F_x)\le2p\qquad\text{for every }x\in{\mathcal S},
\end{equation}
where $F_x$ is the failure event~\eqref{eq:failure-event}: some dilation $t\in[1,2]$ makes all the points $x+tq^n$ with $n\in{\mathcal N}$ miss $B$.

\subsection{From probabilities to one hitting set}

The estimate~\eqref{eq:pointwise-miss} concerns each base point separately, and we cannot take a union bound over uncountably many base points. Instead we integrate the probability that a base point is exceptional, and then cover the remaining base points by a small open set. Several sets appear, with the following roles:
\begin{gridtable}
\begin{tabular}{|>{\raggedright\arraybackslash}p{2.6cm}|>{\raggedright\arraybackslash}p{9.8cm}|}
\hline
\emph{set} & \emph{meaning}\\ \hline
$B$ & the random set constructed from the tables\\ \hline
$B^+$ & an open enlargement of $B$\\ \hline
$R$ & the base points for which some dilation still misses $B^+$ at every index in ${\mathcal N}$\\ \hline
$V$ & an open neighborhood of $R$\\ \hline
$H=B^+\cup V$ & the final hitting set\\ \hline
\end{tabular}
\end{gridtable}
Note that $R$ is a set of base points, not a set of sequence points which missed $B^+$. It also differs from the deterministic set ${\mathbb R}\setminus{\mathcal S}$ of unstable base points: a stable base point can be exceptional for a particular realization of the tables, and an unstable one can nevertheless be handled successfully.

We enlarge $B$ to an open set $B^+$ before defining $R$ for two reasons: the final hitting set $H$ must be open, and the openness of $B^+$ makes its complement closed, which is what allows us to prove that $R$ is closed. For each realization of the tables, enlarge the finitely many selected cells per period to an open periodic set $B^+\supset B$ with
\[
\rho(B^+)\le\rho(B)+p,
\]
which is possible by adding sufficiently short open intervals around the finitely many cell endpoints in each period; we make one fixed such choice for each of the finitely many realizations. Define
\begin{equation}\label{eq:exceptional}
R=\{x:\text{some }t\in[1,2]\text{ has }x+tq^n\notin B^+
\text{ for all }n\in{\mathcal N}\}.
\end{equation}

The set $R$ is periodic, since $B^+$ is. It is also closed. Let $x_j\in R$ with $x_j\longrightarrow x$, and choose dilations $t_j\in[1,2]$ witnessing $x_j\in R$, so that
\[
x_j+t_jq^n\notin B^+\qquad\text{for every }n\in{\mathcal N}.
\]
By compactness of $[1,2]$, pass to a subsequence with $t_j\longrightarrow t\in[1,2]$. For each of the finitely many $n\in{\mathcal N}$, the points $x_j+t_jq^n$ converge to $x+tq^n$ and lie in the closed set ${\mathbb R}\setminus B^+$, so $x+tq^n\notin B^+$. Thus $t$ witnesses $x\in R$. The closedness of $R$ is used only to know that $R$ is a Borel set. This is needed for $\rho(R)$ to be defined, for the function $x\mapsto{\mathbb P}(x\in R)$ to be measurable, so that ${\mathbb E}\rho(R)$ can be computed by integrating it below, and for the outer approximation of $R$ by an open set $V$ at the end of this subsection.

Next we change the direction of averaging: so far we estimated probabilities for each fixed $x$, and now we integrate them over a period. If a stable base point $x$ lies in $R$, then the failure event $F_x$ occurs, since $B\subset B^+$; hence ${\mathbb P}(x\in R)\le{\mathbb P}(F_x)\le2p$ for every $x\in{\mathcal S}$, by~\eqref{eq:pointwise-miss}. The unstable base points have density less than $p$ by~\eqref{eq:unstable-size} and~\eqref{eq:gap-choice}. Hence
\[
{\mathbb E}\rho(R)=\int_0^1{\mathbb P}(x\in R)dx
\le\underbrace{2p\,\rho({\mathcal S})}_{\text{stable base points}}
+\underbrace{\rho({\mathbb R}\setminus{\mathcal S})}_{\text{unstable base points}}
\le3p.
\]
Since there are finitely many realizations of the tables, each giving a closed, hence measurable, set $R$, this interchange is only an interchange of a finite weighted sum and an integral, and no independence between different base points is needed. Moreover ${\mathbb E}\rho(B^+)\le{\mathbb E}\rho(B)+p=2p$ by lemma~\ref{lem:density}.

We now select a good realization for both sets simultaneously. The nonnegative quantity
\[
Y=\rho(B^+)+\rho(R)
\]
satisfies ${\mathbb E}Y\le5p$, so at least one realization has $Y\le5p$. This gives the joint bound needed below; it does not assert that both terms separately satisfy their individual expectation bounds in that same realization. Fix this realization. By outer regularity of Lebesgue measure on the circle, choose an open periodic set $V\supset R$ with $\rho(V)\le\rho(R)+p$, and put
\[
H=B^+\cup V.
\]
It is open and periodic, and
\[
\rho(H)\le\rho(B^+)+\rho(R)+p\le6p.
\]

\subsection{The repair of exceptional base points}

It remains to show that $H$ intersects every normalized copy. The base points outside $R$ are handled by the finite construction, and those in $R$ are handled by the principle already used by Tom in his work with Pramanik~\cite[lemma~0.2]{Tom2015}: an open neighborhood of the base point captures all sufficiently late points of a null sequence. Here we combine it with a quantitative measure allowance.

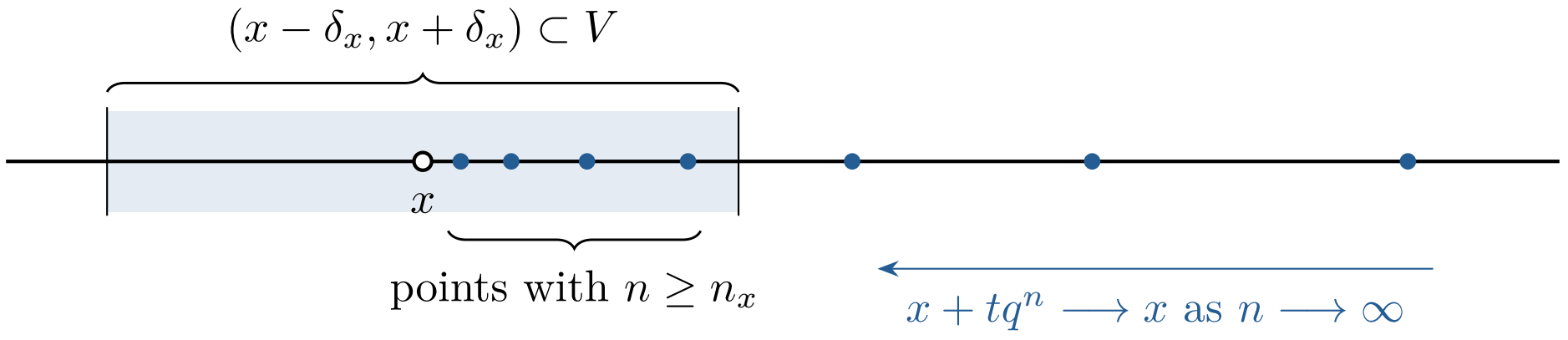
\begin{figure}[htbp]
\centering
\begin{tikzpicture}[x=1cm,y=1cm,>=Stealth,font=\small,
  brace/.style={decorate,decoration={brace,amplitude=4pt},line width=.55pt}]
\fill[leftbranch!12] (.5,-.4) rectangle (5.5,.4);
\draw[thick] (-.3,0)--(12,0);
\draw (.5,-.43)--(.5,.43);
\draw (5.5,-.43)--(5.5,.43);
\draw[black,fill=white,thick] (3.0,0) circle (.07);
\node[below] at (3.0,-.12) {$x$};
\foreach \a in {3.3,3.7,4.3,5.1,6.4,8.3,10.8} {\fill[leftbranch] (\a,0) circle (.065);}
\draw[brace] (.5,.55)--(5.5,.55) node[midway,above=5pt] {$(x-\delta_x,x+\delta_x)\subset V$};
\draw[brace] (5.2,-.55)--(3.2,-.55) node[midway,below=5pt] {points with $n\ge n_x$};
\draw[->,leftbranch] (11,-.85)--(6.6,-.85);
\node[below,text=leftbranch] at (8.8,-.9) {$x+tq^n\longrightarrow x$ as $n\longrightarrow\infty$};
\end{tikzpicture}
\caption{Repairing an exceptional base point $x\in R$. The base point (open circle) is not a point of the affine copy, since $tq^n\ne0$ for every $n$. The points $x+tq^n$ (filled) approach $x$ from the right, and those with $n\ge n_x$ lie in the neighborhood $(x-\delta_x,x+\delta_x)$, hence in $V$. Openness of $V$ is essential: merely adding the set $R$ to $H$ would not justify this conclusion.}
\label{fig:repair}
\end{figure}

Fix $x\in{\mathbb R}$ and $t\in[1,2]$.
\begin{itemize}
\item If $x\notin R$, then by~\eqref{eq:exceptional} some index $n\in{\mathcal N}$ satisfies $x+tq^n\in B^+\subset H$.
\item If $x\in R$, then since $V$ is open, there is $\delta_x>0$ with $(x-\delta_x,x+\delta_x)\subset V$. Choose $n_x$ with $2q^{n_x}<\delta_x$. For every $n\ge n_x$ and every $t\in[1,2]$, we have $0<tq^n\le2q^n<\delta_x$, so $x+tq^n\in V\subset H$; see figure~\ref{fig:repair}.
\end{itemize}
In the second case the repair works for all normalized dilations at once, and the indices $n\ge n_x$ need not belong to ${\mathcal N}$, which does not matter, since any natural number $n$ is permitted. In either case the copy intersects $H$. This proves proposition~\ref{prop:periodic}, and by section~\ref{sec:periodic} it completes the proof of theorem~\ref{thm:main}: the normalized hitting sets are combined there by the summable construction~\eqref{eq:global-budget}, following the organization by periodic hitting sets shared with the Rajchman paper~\cite[proposition~3.2]{IKMRY2026Rajchman}.

\subsection{Choosing a realization by conditional expectations}

The choice of a good realization can also be made deterministically, by the classical method of conditional expectations. Alon and Spencer~\cite[section~16.1]{AlonSpencer2008} present this method and trace its development through the work of Erd\H{o}s and Selfridge, Spencer, and Raghavan. Written as a nonincreasing potential, it is a second instance of descent in this proof.

Fix the enlargements $B^+$ chosen above and the cost $Y=\rho(B^+)+\rho(R)$, so that ${\mathbb E}Y\le5p$, and list every random table entry in some order. After the first $\ell$ entries have been chosen, finitely many completions remain, and we let $\Phi_\ell$ be the conditional expectation of $Y$ given those choices, that is, the weighted average of the cost over the remaining completions. If the next entry is one with probability $a$, its two possible continuations have conditional expectations $z_0$ and $z_1$ satisfying
\[
\Phi_\ell=(1-a)z_0+az_1.
\]
At least one continuation is therefore at most $\Phi_\ell$; choose it. Repeating this choice produces
\[
0\le\Phi_{\ell+1}\le\Phi_\ell\le\Phi_0\le5p,
\]
and after the final entry has been chosen, the conditional expectation is the actual value of $Y$.

The two potentials in the proof should be distinguished. Section~\ref{sec:descent} tracks a routing probability as the depth increases, while $\Phi_\ell$ tracks the expected final measure cost as table entries are fixed. The latter gives a deterministic selection rule in principle, although evaluating these conditional expectations may be very expensive; the random construction and its local geometric estimates remain essential to the initial bound $\Phi_0\le5p$.

\section{A closer comparison with the two earlier papers}\label{sec:comparison}

We now return to the two-fold sumset paper and the Rajchman paper, both cited in the OpenAI manuscript, and describe what each of them achieves.
\begin{itemize}
\item The two-fold sumset paper of Mora Cu\'ellar, Iosevich, Kulkarni, Rojas Aravena, and Yavicoli~\cite{IKMRY2026Sumset} proves that if $A\subset{\mathbb R}$ is any infinite set, $a\ne0$ and $0<|\beta|<1$, then neither $\{a\beta^n:n\ge1\}+A$ nor $\{a\beta^n:n\ge1\}-A$ is measure universal; more generally, the geometric sequence may be replaced by any set containing a lacunary sequence $(b_n)_{n\in{\mathbb N}}$ with $-\log b_n=O(n)$; here a positive sequence is lacunary if $b_{n+1}\le\lambda b_n$ for some fixed $\lambda\in(0,1)$. Its method uses additive structure to produce many separated points.
\item The Rajchman paper of Iosevich, Kulkarni, Mora Cu\'ellar, Rojas Aravena, and Yavicoli~\cite{IKMRY2026Rajchman} proves that if a set supports a probability measure whose Fourier transform tends to zero at infinity, then for every $\varepsilon\in(0,1)$ there is a closed, $1$-periodic, nowhere dense set, of measure at least $1-\varepsilon$ in every interval of length one, containing no affine copy of it. Its method uses Fourier decay to construct small periodic sets which intersect every affine copy.
\end{itemize}
This section examines these two mechanisms and their relationship with the geometric construction. The connections are most informative when stated at the level of individual estimates: a similar organization does not by itself show that one theorem follows from another, nor does it establish how the later proof was discovered.

\subsection{Separated points and the cost of discretization}\label{subsec:separated}

Recall the quantities of section~\ref{sec:obstacle}. For a finite set of $N$ points, let
\[
\delta=\frac{\text{smallest distance between distinct points}}{\text{diameter of the finite set}}.
\]
Separation has two consequences, which pull in opposite directions.
\begin{itemize}
\item Separation allows the $N$ points to occupy different random cells, so that at one fixed dilation they all miss the selected cells with probability $(1-p)^N\le e^{-pN}$.
\item Finer separation requires finer cells, which increases the number of cell configurations encountered as the dilation varies.
\end{itemize}
Combining them gives a bound of the form
\[
C\,\frac{N}{\delta}\,e^{-pN}
\]
for the probability that some dilation fails, and taking logarithms makes the competition explicit:
\[
\underbrace{\log C+\log N+\log\frac{1}{\delta}}_{\text{cost of allowing all cell configurations}}
\;-\;\underbrace{pN}_{\text{gain from independent tests}}.
\]
The condition
\[
\log\frac{1}{\delta}=o(N)
\]
now has a clear meaning: the required resolution grows slowly enough, on a logarithmic scale, for the number of available points to outweigh it, even when $p$ is small. If arbitrarily large finite subsets satisfy this condition, the set is not measure universal. This criterion for finite subsets originates with Kolountzakis~\cite{Kolountzakis1997}, and its translation-invariant formulation is due to Chleb\'ik~\cite[theorem~15]{Chlebik2015}; the two-fold sumset paper builds on it~\cite[theorem~9 and proposition~14]{IKMRY2026Sumset}. The probabilistic calculation behind it can be seen directly in Chleb\'ik's proof: for a fixed base point, at most a constant times $\frac{N}{\delta}$ dilation representatives are needed, and the chance of missing all $N$ separated points at one representative is $(1-p)^N$.

For a single geometric sequence the criterion fails, by~\eqref{eq:gap-obstruction}. The sumset paper produces the required excess of points from the additional summand, by controlling how often two cross-sums can be close. Given two positive sequences $(b_i)_{i\in{\mathbb N}}$ and $(d_j)_{j\in{\mathbb N}}$, consider the table of cross-sums
\[
s_{ij}=b_i+d_j.
\]
The inequality
\[
|(b_i-b_{i'})+(d_j-d_{j'})|\le e^{-W}
\qquad\text{simply says}\qquad
|s_{ij}-s_{i'j'}|\le e^{-W},
\]
so counting the quadruples $(i,i',j,j')$ which satisfy it counts the pairs of entries of the table which coincide or lie very close together. This is the meaning of the near-additive-energy estimate. Recall that a positive sequence $(b_i)_{i\in{\mathbb N}}$ is lacunary if $b_{i+1}\le\lambda b_i$ for some fixed $\lambda\in(0,1)$. After passing to subsequences whose successive ratios are below $\frac{1}{2}$, lemma~19 of the sumset paper bounds the number of such quadruples by a constant times $I(W)J(W)$, where $I(W)$ and $J(W)$ count the terms of the respective sequences which are at least $e^{-W}$, and all four indices range over those terms.

A worked example shows why the extra summand helps. Take $b_i=d_i=4^{-i}$ for $1\le i\le\ell$. For $\ell=3$ the table of cross-sums, multiplied by $64$, is
\begin{gridtable}
\begin{tabular}{|c|c|c|c|}
\hline
$64\,s_{ij}$ & $d_1=\tfrac{1}{4}$ & $d_2=\tfrac{1}{16}$ & $d_3=\tfrac{1}{64}$\\ \hline
$b_1=\tfrac{1}{4}$ & $32$ & $20$ & $17$\\ \hline
$b_2=\tfrac{1}{16}$ & $20$ & $8$ & $5$\\ \hline
$b_3=\tfrac{1}{64}$ & $17$ & $5$ & $2$\\ \hline
\end{tabular}
\end{gridtable}
There are nine indexed pairs but only six distinct sums; the repeated entries show why counting pairs alone is not enough, and why a bound on near-coincidences is needed. For general $\ell$, base-four expansions show that the unordered pairs give exactly $N=\frac{\ell(\ell+1)}{2}$ distinct sums. These sums are multiples of $4^{-\ell}$, so their minimum separation is at least $4^{-\ell}$, while their diameter is bounded above and below independently of $\ell$. Consequently
\[
\log\frac{1}{\delta}\le\ell\log4+O(1),\qquad N\asymp\ell^2,
\qquad\text{so}\qquad\log\frac{1}{\delta}=o(N).
\]
The number of distinct points grows quadratically, while the logarithmic resolution cost grows only linearly. This simple case illustrates the mechanism; the cited estimates handle far more general sequences.

In general, a bound on near-coincidences yields many separated sums. One way to see this is a binning argument. Divide the line into intervals of length $e^{-W}$, and let $n_k$ be the number of indexed sums in the $k$th interval, so that $\sum_kn_k=I(W)J(W)$. Two sums in the same interval are within $e^{-W}$ of each other, so $\sum_kn_k^2$ is at most the number of near-coincident pairs, which is at most $C\,I(W)J(W)$. By the Cauchy--Schwarz inequality,
\[
\bigl(I(W)J(W)\bigr)^2=\Bigl(\sum_kn_k\Bigr)^2\le\#\{k:n_k>0\}\sum_kn_k^2,
\]
so at least $I(W)J(W)/C$ intervals are occupied. Selecting one sum from each occupied interval of one parity gives at least $I(W)J(W)/(2C)$ sums whose mutual distances are at least $e^{-W}$. The proof of proposition~20 of the sumset paper extracts a comparable number of separated cross-sums, and the hypothesis $\limsup_{W\to\infty}\frac{I(W)J(W)}{W}=\infty$ of that proposition then supplies the criterion for finite subsets. Indeed, the selected sums have mutual distances at least $e^{-W}$ and a diameter bounded independently of $W$, so that $\log\frac{1}{\delta}\le W+O(1)$, while along a suitable sequence of values of $W$ the number of points is a large multiple of $W$.

The geometric proof faces the same comparison between the number of tests and the logarithm of the number of cell configurations met as the dilation varies. The two estimates can be placed side by side; the last row records upper bounds, which is all that the union bound requires:
\begin{gridtable}
\begin{tabular}{|>{\raggedright\arraybackslash}p{4.2cm}|>{\raggedright\arraybackslash}p{3.6cm}|>{\raggedright\arraybackslash}p{4.6cm}|}
\hline
& \emph{separated-point argument} & \emph{tree argument}\\ \hline
number of tests & $N$ & $(M-1)r$\\ \hline
failure bound at one dilation & $e^{-pN}$ & $e^{-p(M-1)r/2}$\\ \hline
upper bound for the logarithm of the number of cell configurations & $\log N+\log\frac{1}{\delta}+O(1)$ & $2r\log\frac{1}{q}+\log r+O_{M,q}(1)$\\ \hline
\end{tabular}
\end{gridtable}
In the tree argument $M$ is fixed before $r$ increases. The significant feature is that increasing $M$ improves the coefficient of $r$ in the failure exponent without increasing the coefficient $2\log\frac{1}{q}$ of the leading term in the upper bound for the logarithm of the number of cell configurations; the bound on the span of each child subtree is what keeps $M$ out of that coefficient. This is a concrete mathematical parallel: in one proof additive structure supplies many separated cross-sums, while in the other routing supplies many tests whose cell configurations can still be counted locally.

Finally, the relationship between the two theorems should be stated accurately. For any $a_0\in A$,
\[
a_0+G_q\subseteq G_q+A,
\]
so every affine copy of $G_q+A$ contains an affine copy of $G_q$. A set avoiding all affine copies of $G_q$ therefore also avoids all affine copies of $G_q+A$, but avoiding the larger pattern alone does not imply avoiding the smaller one. In particular, theorem~\ref{thm:main} now implies the main statement of the sumset paper~\cite{IKMRY2026Sumset}, which concerns the sets $\{a\beta^n:n\ge1\}\pm A$ with $a\ne0$ and $0<|\beta|<1$. If $0<\beta<1$, these sets contain affine copies of $G_\beta$. If $-1<\beta<0$, the even powers satisfy $\beta^{2n}=(\beta^2)^n$ with $0<\beta^2<1$, so these sets contain affine copies of $G_{\beta^2}$. The sumset paper keeps its significance for other reasons. Its main theorem was proved earlier, independently of the tree construction and by a different method, and its near-additive-energy estimate is not an input to the tree proof. Its broader results also concern sequences which are not geometric, about which theorem~\ref{thm:main} says nothing unless the set happens to contain an affine copy of some $G_q$: summands which only contain a lacunary sequence with $-\log b_n=O(n)$, sums and differences $S_1\pm S_2$ of sets with lacunary subsequences whose counting functions satisfy $\limsup_{W\to\infty}\frac{I(W)J(W)}{W}=\infty$, and examples such as $\{2^{-n^{\alpha_1}}\}+\{2^{-n^{\alpha_2}}\}$ with $\frac{1}{\alpha_1}+\frac{1}{\alpha_2}>1$.

\subsection{Periodic hitting sets and the Rajchman argument}\label{subsec:rajchman}

A finite positive Borel measure $\mu$ on ${\mathbb R}$ is called a Rajchman measure if its Fourier transform
\[
\widehat\mu(\zeta)=\int_{{\mathbb R}}e^{-2\pi i\zeta u}d\mu(u)
\]
tends to zero as $|\zeta|\longrightarrow\infty$. The Rajchman paper turns this Fourier decay into small hitting sets through the following chain of implications:
\[
\text{Fourier decay}\ \Longrightarrow\ \text{uniform distribution of dilates modulo one}\ \Longrightarrow\ \text{small periodic hitting sets}.
\]

\subsubsection*{From Fourier decay to uniform distribution}

``Modulo one'' means working on the circle ${\mathbb T}={\mathbb R}/{\mathbb Z}$, the unit interval with its endpoints identified. A probability measure $\mu$ supported on a compact set $K$ can be viewed as the distribution of a random point $u\in K$, and for $\tau>0$ we consider the distribution of $\tau u$ modulo one. The relevant conclusion is that, for every continuous $1$-periodic function $f$,
\begin{equation}\label{eq:equidistribution}
\int f(\tau u)\,d\mu(u)\longrightarrow\int_0^1f(v)\,dv\qquad\text{when }\tau\longrightarrow\infty.
\end{equation}
This is what convergence to the uniform measure on the circle means; it does not require $\mu$ to have a density. The proof has three steps.
\begin{itemize}
\item For $f(v)=e^{-2\pi ijv}$ with an integer $j$, the left side of~\eqref{eq:equidistribution} equals $\widehat\mu(j\tau)$.
\item For $j\ne0$, Fourier decay makes $\widehat\mu(j\tau)$ tend to zero, which is the value of the uniform integral; for $j=0$, both sides equal one.
\item Approximation of continuous $1$-periodic functions by trigonometric polynomials, which is classical (Fej\'er's theorem), extends~\eqref{eq:equidistribution} to all continuous $f$.
\end{itemize}
In detail, given $\varepsilon>0$ choose a trigonometric polynomial $P(v)=\sum_{|j|\le J}c_je^{-2\pi ijv}$ with $\|f-P\|_\infty<\varepsilon$. Since $\int P(\tau u)\,d\mu(u)=\sum_{|j|\le J}c_j\widehat\mu(j\tau)$, $\widehat\mu(0)=1$ and $\int_0^1P(v)\,dv=c_0$,
\[
\Bigl|\int f(\tau u)\,d\mu(u)-\int_0^1f(v)\,dv\Bigr|\le2\varepsilon+\sum_{0<|j|\le J}|c_j|\,|\widehat\mu(j\tau)|.
\]
The finite sum tends to zero as $\tau\longrightarrow\infty$, and then $\varepsilon$ may be taken arbitrarily small.
Convergence of measures then gives density of their supports. Suppose that for some $\delta>0$ there were arbitrarily large $\tau$ for which $\tau K$ misses an arc of length $2\delta$ modulo one. Passing to a subsequence along which the centers of these arcs converge, all these dilates eventually miss a fixed arc of length $\delta$. A nonnegative continuous function supported inside that arc and not identically zero has integral zero against the dilated measures but a positive integral against the uniform measure, contradicting~\eqref{eq:equidistribution}. This is lemma~2.3 of the Rajchman paper~\cite{IKMRY2026Rajchman}: if a compact set $K$ supports a Rajchman probability measure, then for every $\delta>0$ all sufficiently large dilates $\tau K$ are $\delta$-dense modulo one, meaning that every point of the circle is within distance $\delta$ of some point of $\tau K$.

\subsubsection*{From uniform distribution to periodic hitting sets}

\begin{figure}[htbp]
\centering
\begin{tikzpicture}[x=1cm,y=1cm,>=Stealth,font=\small]
\draw[thick] (-.3,0)--(12.3,0);
\foreach \m in {0,...,5}{
  \fill[leftbranch!25] ({2.4*\m-.36},-.22) rectangle ({2.4*\m+.36},.22);
  \draw (2.4*\m,-.1)--(2.4*\m,.1);
  \ifnum\m=0 \node[below,font=\footnotesize] at (2.4*\m,-.25) {$0$};
  \else\ifnum\m=5 \node[below,font=\footnotesize] at (2.4*\m,-.25) {$1$};
  \else \node[below,font=\footnotesize] at (2.4*\m,-.25) {$\frac{\m}{5}$};\fi\fi
}
\draw[|-|,black!70] (9.6,.45)--(9.96,.45);
\node[above,font=\footnotesize] at (9.78,.5) {$\frac{\delta}{Q}$};
\fill (1.2,0) circle (.06);
\node[below] at (1.2,-.12) {$x$};
\fill[rightbranch] (7.05,0) circle (.07);
\draw[->,rightbranch,thick] (1.2,.35) to[out=30,in=150] node[midway,above,font=\footnotesize] {translate by $sk$} (7.0,.35);
\node[below,font=\footnotesize,text=rightbranch] at (7.05,-.8) {$x+sk$};
\node[anchor=west,font=\footnotesize] at (-.3,-1.55) {$B_{Q,\delta}$ with $Q=5$: $Q$ intervals per period, each of length $\frac{2\delta}{Q}$, total density $2\delta$};
\end{tikzpicture}
\caption{The periodic set $B_{Q,\delta}$ of the Rajchman paper, shown on one period with $Q=5$. Choosing $k\in K$ so that $Qsk$ is close to $-Qx$ modulo one makes $Q(x+sk)$ close to an integer $m$, that is, $x+sk$ close to $\frac{m}{Q}$, so the translated point lands in one of the shaded intervals.}
\label{fig:periodic-hitting}
\end{figure}
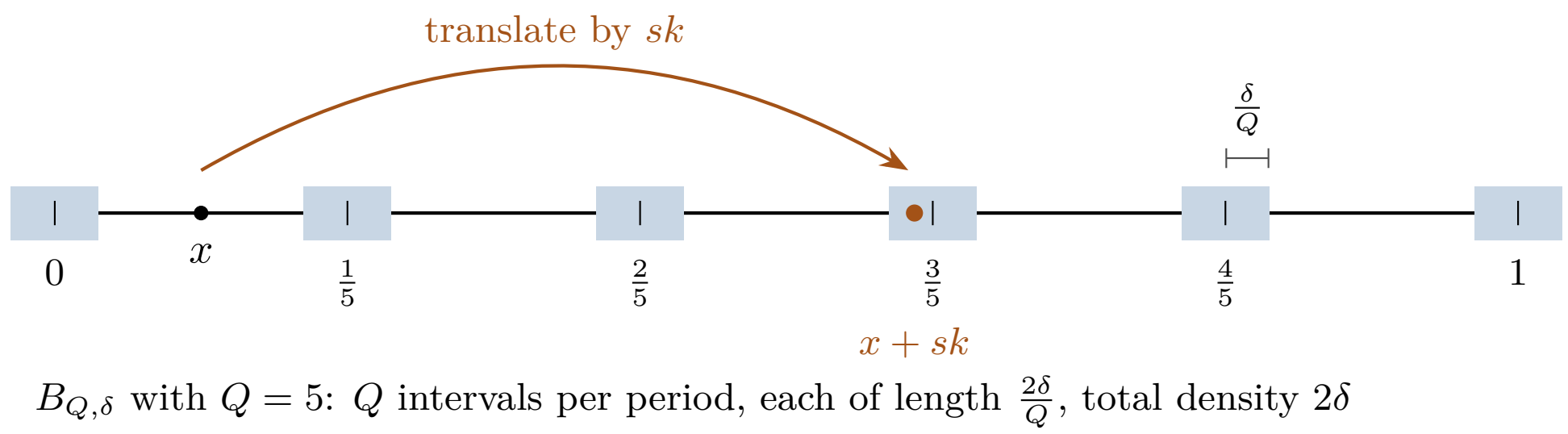

To avoid confusing an integer frequency with our geometric ratio $q$, write the periodic sets of the Rajchman paper as
\[
B_{Q,\delta}=\{u:\operatorname{dist}(Qu,{\mathbb Z})<\delta\}
=\bigcup_{m\in{\mathbb Z}}\Bigl(\frac{m}{Q}-\frac{\delta}{Q},\frac{m}{Q}+\frac{\delta}{Q}\Bigr),
\]
where $Q\ge1$ is an integer and $0<\delta<\frac{1}{2}$. The set consists of narrow intervals around the points $\frac{m}{Q}$. There are $Q$ intervals per unit length, each of length $\frac{2\delta}{Q}$, so its density is $Q\cdot\frac{2\delta}{Q}=2\delta$. Increasing $Q$ makes the intervals more frequent and narrower while preserving their total density; see figure~\ref{fig:periodic-hitting}.

Lemma~3.1 of the Rajchman paper converts density of dilates into a hitting statement. First take $s>0$, and suppose that $sQK$ is $\frac{\delta}{2}$-dense modulo one. Then some $k\in K$ has $sQk$ within $\frac{\delta}{2}$ of the target $-Qx$ modulo one, so that
\[
\operatorname{dist}(Q(x+sk),{\mathbb Z})\le\frac{\delta}{2}<\delta,
\]
that is, $x+sk\in B_{Q,\delta}$. The target depends on $x$, but the density assertion works for every target, hence for every translation. Negative $s$ are handled by reflection, using the target $Qx$ and the density of $|s|QK$.

\subsubsection*{Assembling the dilation ranges}

The proof of proposition~3.2 of the Rajchman paper assigns small summable widths to a sequence of dilation ranges, and chooses integer frequencies large enough that the density of dilates holds on each range. Every nonzero dilation is then covered, while the total removed measure stays small. The same two tasks appear in sections~\ref{sec:periodic} and~\ref{sec:finish} above, where we first construct a small periodic hitting set for the normalized dilations and then combine dilation ranges by means of the summable bound~\eqref{eq:global-budget}. What has changed is the input to the local construction, which is now the estimate~\eqref{eq:all-scales} together with the repair of exceptional base points, rather than Fourier equidistribution.

\subsection{The Fourier hypothesis and geometric sequences}

Two separate facts explain why the Rajchman theorem does not already handle the geometric sequence.
\begin{enumerate}
\item A countable set cannot support a Rajchman probability measure.
\item For the dyadic sequence, even the weaker hypothesis that all sufficiently large dilates are dense modulo one fails.
\end{enumerate}
The first fact does not by itself imply the second. As remark~3.3 of the Rajchman paper explains, the density assumption used there is weaker than the existence of a Rajchman measure, so countability alone does not rule out every use of its construction of periodic hitting sets.

The first fact is a classical consequence of Wiener's theorem: every Rajchman probability measure is atomless, meaning that every singleton has measure zero. The Rajchman paper recalls this as proposition~4.1~\cite{IKMRY2026Rajchman}, following Wiener's argument, and we include the short proof. The guiding idea is that averaging the Fourier transform against the appropriate oscillation detects the mass at a single point. Suppose that $\mu(\{a\})>0$ for some $a$. By Fubini's theorem,
\[
\frac{1}{2T}\int_{-T}^{T}e^{2\pi i\zeta a}\widehat\mu(\zeta)d\zeta
=\int_{{\mathbb R}}\Bigl(\frac{1}{2T}\int_{-T}^{T}e^{-2\pi i\zeta(u-a)}d\zeta\Bigr)d\mu(u),
\]
and the inner average is
\[
\frac{1}{2T}\int_{-T}^{T}e^{-2\pi i\zeta(u-a)}d\zeta=
\begin{cases}
1, & u=a,\\[3pt]
\dfrac{\sin(2\pi T(u-a))}{2\pi T(u-a)}, & u\ne a.
\end{cases}
\]
It is bounded by one and converges pointwise to the indicator function of $\{a\}$ as $T\longrightarrow\infty$, so by dominated convergence the right side tends to $\mu(\{a\})$. The Rajchman condition makes the left side tend to zero: since $|\widehat\mu|\le1$, for $0<R_0<T$,
\[
\Bigl|\frac{1}{2T}\int_{-T}^{T}e^{2\pi i\zeta a}\widehat\mu(\zeta)\,d\zeta\Bigr|
\le\frac{R_0}{T}+\sup_{|\zeta|\ge R_0}|\widehat\mu(\zeta)|,
\]
and letting first $T\longrightarrow\infty$ and then $R_0\longrightarrow\infty$ gives the limit zero. Thus $\mu(\{a\})=0$ for every $a$. If $K$ were countable and supported such a probability measure, we would have
\[
1=\mu(K)=\sum_{a\in K}\mu(\{a\})=0,
\]
which is impossible. In particular, neither $G_q$ nor $G_q\cup\{0\}$ supports a Rajchman probability measure.

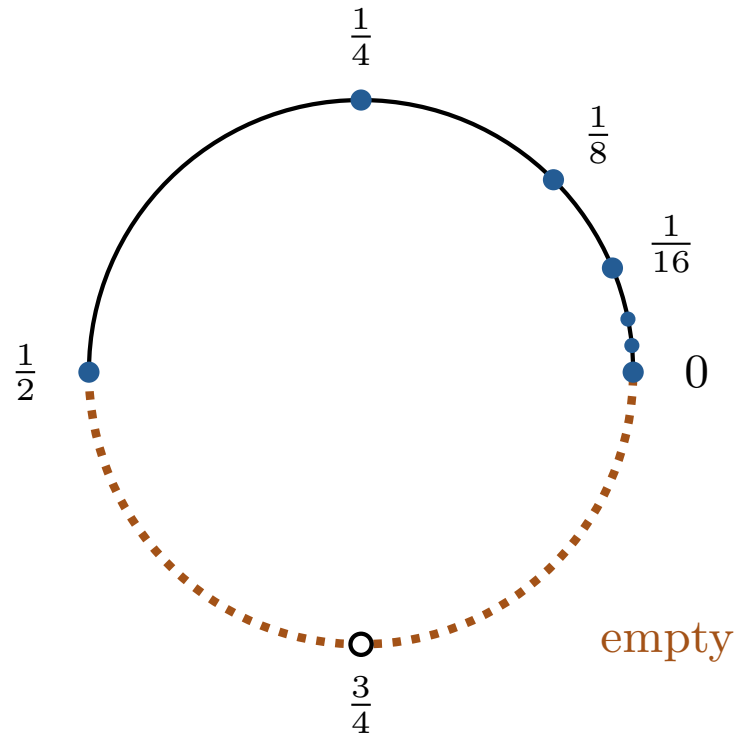
\begin{figure}[htbp]
\centering
\begin{tikzpicture}[font=\small,>=Stealth]
\def\R{1.8}
\draw[thick] (\R,0) arc (0:180:\R);
\draw[rightbranch,line width=1.6pt,dotted] (-\R,0) arc (180:360:\R);
\foreach \v/\lab in {0/0,0.5/\frac{1}{2},0.25/\frac{1}{4},0.125/\frac{1}{8},0.0625/\frac{1}{16}}{
  \fill[leftbranch] ({360*\v}:\R) circle (.07);
  \node[font=\footnotesize] at ({360*\v}:{\R+.42}) {$\lab$};
}
\foreach \v in {0.03125,0.015625}{\fill[leftbranch] ({360*\v}:\R) circle (.05);}
\draw[black,fill=white,thick] (270:\R) circle (.07);
\node[font=\footnotesize] at (270:{\R+.4}) {$\frac{3}{4}$};
\node[align=left,anchor=north west,font=\footnotesize,text=rightbranch] at (315:{\R+.25}) {empty arc $(\frac{1}{2},1)$};
\node[anchor=west,align=left,font=\footnotesize] at (2.8,.3) {$K=\{0,\frac{1}{2},\frac{1}{4},\frac{1}{8},\ldots\}$ on ${\mathbb T}$;\\the picture is unchanged\\after multiplying by $2^j$};
\end{tikzpicture}
\caption{The dyadic set $K$ on the circle ${\mathbb R}/{\mathbb Z}$, with $0$ on the right and positions measured counterclockwise. Every dilate $2^jK$ coincides with $K$ modulo one, so no dilate intersects the empty arc $(\frac{1}{2},1)$, drawn dotted, and the point $\frac{3}{4}$ stays at distance $\frac{1}{4}$ from all of them.}
\label{fig:dyadic-circle}
\end{figure}

The second fact is shown by the basic dyadic example
\[
K=\{0\}\cup\{2^{-n}:n\ge1\}=\Bigl\{0,\frac{1}{2},\frac{1}{4},\frac{1}{8},\ldots\Bigr\}.
\]
Multiplication by $8$ sends the first three nonzero terms to the integers $4$, $2$, $1$, while
\[
8\cdot\frac{1}{16}=\frac{1}{2},\qquad8\cdot\frac{1}{32}=\frac{1}{4},\qquad\ldots
\]
reproduces the sequence. Thus $8K=K$ modulo one, and in the same way $2^jK=K$ modulo one for every positive integer $j$. The point $\frac{3}{4}$ is always at circular distance $\frac{1}{4}$ from this set, so these arbitrarily large dilates miss the arc $(\frac{1}{2},1)$ and fail, for example, to be $\frac{1}{8}$-dense; see figure~\ref{fig:dyadic-circle}. This shows that the sufficient density hypothesis fails for the dyadic example; it does not rule out other periodic hitting constructions. The tree construction requires no such equidistribution.

There is also a difference in the final global sets. The Rajchman construction uses integer frequencies, so all its hitting sets remain $1$-periodic, and its avoiding set is $1$-periodic, as in its main theorem. Our normalized hitting sets are $1$-periodic, but the union of their signed dilations in section~\ref{sec:periodic} need not be, and we have used exactly the measure estimate needed on $[0,1]$.

\FloatBarrier
\subsection{What the tree contributes}

All three approaches share one strategy: turn avoidance into the problem of intersecting every copy with a small set, and combine the required dilation ranges with a controlled total measure. In each of them the heart of the matter is what we call the \emph{local hitting property}: for one bounded range of dilations, a periodic set of small density intersects every copy whose dilation lies in that range. For the tree construction this is proposition~\ref{prop:periodic}, for the dilations $t\in[1,2]$, and the dilation ranges are then combined as in section~\ref{sec:periodic}. The random-cell arguments, including the sumset construction, and the tree construction share a second feature: both compare the number of useful tests with the number of cell configurations met as the dilation varies, as in the table at the end of section~\ref{subsec:separated}. The three approaches differ in how they obtain the local hitting property:
\begin{gridtable}
\begin{tabular}{|>{\raggedright\arraybackslash}p{3.4cm}|>{\raggedright\arraybackslash}p{9.2cm}|}
\hline
\emph{approach} & \emph{how the local hitting property is obtained}\\ \hline
random cells and separated points~\cite{Kolountzakis1997,Chlebik2015,IKMRY2026Sumset}; section~\ref{subsec:separated} & Periodic cells are selected independently with probability $p$. If $N$ points of a copy lie in distinct cells, they all miss the selected cells with probability at most $e^{-pN}$. A union bound over the cell configurations met as the dilation varies then succeeds when $N$ is large compared with the logarithm of the number of these configurations. In the sumset paper, the second summand supplies the required separated points.\\ \hline
Fourier decay~\cite{IKMRY2026Rajchman}; section~\ref{subsec:rajchman} & No randomness is used. Fourier decay of a measure supported on the pattern $K$ makes all the large dilates of $K$ $\delta$-dense modulo one. Consequently, for all dilations $s$ in a suitable range, the periodic set $B_{Q,\delta}$, formed by narrow intervals around the points $\frac{m}{Q}$, intersects every translate $x+sK$.\\ \hline
random tree~\cite{OpenAI2026}; sections~\ref{sec:periodic}--\ref{sec:finish} & Fair random choices route each point down a finite tree, and a final random choice at the leaf decides, with probability $p$, whether the point belongs to the random set. Its expected density of the random set is $p$, however large the tree is (lemma~\ref{lem:density}), while the branching gives $(M-1)r$ local tests, which all fail at one dilation with probability $(1-\frac{p}{2})^{(M-1)r}$ (lemma~\ref{lem:failure}). Each test depends only on its own edge and the subtree below it, so at most $C_{M,q}\,r\,q^{-2r}$ representative dilations need to be checked (section~\ref{sec:scales}). The union bound succeeds once $\frac{p(M-1)}{2}>2\log\frac{1}{q}$, and the few exceptional base points are covered by a small open set (section~\ref{sec:finish}).\\ \hline
\end{tabular}
\end{gridtable}
What makes the tree construction work is the combination of three facts proved above. First, a point belongs to the random set with probability $p$, whatever the size of the tree (lemma~\ref{lem:density}). Second, at a default vertex the nearby points of the copy supply $(M-1)r$ local tests (section~\ref{sec:routing} and lemma~\ref{lem:failure}). Third, these tests change with the dilation only at the crossing times of their local grids; since the group of each edge $e_i$ spans at most $2r$ indices (lemma~\ref{lem:span}), each tested point has at most $1+2cq^{-2r}$ crossing times (section~\ref{subsec:crossings}), and hence at most $D_{M,q}(r)$ representative dilations need to be checked.

The inequality~\eqref{eq:alpha}, $\frac{p(M-1)}{2}>2\log\frac{1}{q}$, summarizes the resulting advantage. Its left side is the rate, in $r$, at which the probability that all the local tests fail at one dilation decreases, and its right side is the rate at which the number of representative dilations grows, by~\eqref{eq:rate}. The other parameters each have one task: the depth $d$ makes a default likely, by~\eqref{eq:depth-choice}; the gap $g$ makes the earlier routing decisions stable, by lemma~\ref{lem:stable} and~\eqref{eq:gap-choice}; and the lower bound $r_0$ on the interval lengths makes the remaining failure probability small, by~\eqref{eq:window-choice}. The order in which these parameters are chosen in section~\ref{sec:finish} ensures that their tasks are compatible.

\end{document}